\input amstex
\magnification=1200
\documentstyle{amsppt}
\vsize 18.8cm

\def\dist{\operatorname {\roman {dist}}}

\def\eps{\epsilon}
\def\weakto{\rightharpoonup}
\def\om{\Omega}
\def\co{\colon}
\def\R{\Bbb R}
\def\d{\roman{d}}

\def \R{\Bbb R}
\def\cl{\bar}
\def\N{\Bbb N}

\newcount\fornumber\newcount\artnumber\newcount\tnumber
\newcount\secnumber
\def\Ref#1{%
  \expandafter\ifx\csname mcw#1\endcsname \relax
    \warning{\string\Ref\string{\string#1\string}?}%
    \hbox{$???$}%
  \else \csname mcw#1\endcsname \fi}
\def\Refpage#1{%
  \expandafter\ifx\csname dw#1\endcsname \relax
    \warning{\string\Refpage\string{\string#1\string}?}%
    \hbox{$???$}%
  \else \csname dw#1\endcsname \fi}
\def\warning#1{\immediate\write16{
            -- warning -- #1}}
\def\CrossWord#1#2#3{%
  \def\x{}%
  \def\y{#2}%
  \ifx \x\y \def\z{#3}\else
            \def\z{#2}\fi
  \expandafter\edef\csname mcw#1\endcsname{\z}
\expandafter\edef\csname dw#1\endcsname{#3}}
\def\Tag#1#2{\begingroup
  \edef\mmhead{\string\CrossWord{#1}{#2}}%
  \def\writeref{\write\refout}%
  \expandafter \expandafter \expandafter
  \writeref\expandafter{\mmhead{\the\pageno}}%
\endgroup}
\openin15=\jobname.ref
\ifeof15 \immediate\write16{No file \jobname.ref}%
\else                   \input \jobname.ref \fi \closein15
\newwrite\refout
\openout\refout=\jobname.ref
\def\newsection{\global\advance\secnumber by 1\fornumber=0\tnumber=0}
\def\dff#1..{%
     \global\advance\fornumber by 1
     \Tag{#1}{(\the\secnumber.\the\fornumber)}
 \unskip                   \the\secnumber.\the\fornumber
   \ignorespaces
}%
\def\dffs#1..{%
     \global\advance\fornumber by 1
     \Tag{#1}{\the\secnumber.\the\fornumber}
 \unskip                   \the\secnumber.\the\fornumber
   \ignorespaces
}%
\def\rff#1..{\Ref{#1}}
\def\dft#1..{%
     \global\advance\tnumber by 1
     \Tag{#1}{\the\secnumber.\the\tnumber}
 \unskip                   \the\secnumber.\the\tnumber
   \ignorespaces
}%
\def\rft#1..{\Ref{#1}}
\def\dfa#1..{%
     \global\advance\artnumber by 1
     \Tag{#1}{\the\artnumber}
 \unskip                   \the\artnumber
   \ignorespaces
}%
\def\rfa#1..{\Ref{#1}}
\def\dep#1..{\Tag{#1}{}}
\def\rep#1..{p. \Refpage{#1}}

\topmatter

\title
Attractors and global averaging of non-autonomous reaction-diffusion equations
in $\R^N$
\endtitle

\author
Francesca Antoci\\and\\ Martino Prizzi
\endauthor

\leftheadtext{Francesca Antoci and Martino Prizzi}
\rightheadtext{Attractors and global averaging}

\address
Francesca Antoci, Politecnico di Torino, Dipartimento di Matematica,
corso Duca degli Abruzzi 24, 10129 Torino, Italy
\endaddress

\email
antoci\@calvino.polito.it
\endemail

\address
Martino Prizzi, Universit\`a degli Studi di Trieste, Dipartimento
di Scienze Matematiche, via Valerio 12, 34127 Trieste, Italy
\endaddress

\email
prizzi\@mathsun1.univ.trieste.it
\endemail

\date
\enddate

\abstract We consider a family of non-autonomous
reaction-diffusion equations $$ u_t=\sum_{i,j=1}^N a_{ij}(\omega
t)\partial_i\partial_j u+f(\omega t,u)+ g(\omega t,x),\quad (t,
x)\in\R_+\times\R^N\tag{$E_\omega$} $$ with almost periodic,
rapidly oscillating principal part and nonlinear interactions. As
$\omega\to +\infty$, we prove that the solutions of $(E_\omega)$
converge to the solutions of the averaged equation $$
u_t=\sum_{i,j=1}^N \cl a_{ij}\partial_i\partial_j u+\cl f(u)+ \cl
g(x),\quad(t, x)\in\R_+\times\R^N.\tag{$E_\infty$} $$ If $f$ is
dissipative, we prove existence and upper-semicontinuity of
attractors for the family $(E_\omega)$ as $\omega\to+\infty$.
\endabstract

\endtopmatter

\document

\newsection\head \the\secnumber. Introduction\endhead
In this paper we study a family of non-autonomous
reaction-diffusion equations $$ u_t=\sum_{i,j=1}^N a_{ij}(\omega
t)\partial_i\partial_j u+f(\omega t,u)+ g(\omega t,x),\quad (t,
x)\in\R_+\times\R^N\tag\dff intro1.. $$ with almost periodic,
rapidly oscillating principal part and nonlinear interactions.
Under suitable hypothesis (see Section 2), the Cauchy problem for
\rff intro1.. is well-posed in $H^1(\R^N)$ and the equation
generates a (global) process, that is, a 
two-parameter family of nonlinear operators $\Pi_\omega(t,s)$ from $H^1(\R^N)$ into itself 
such that $$
\cases \Pi_\omega(t,p)\Pi_\omega(p,s)=\Pi_\omega(t,s)&t\geq p \geq
s\\ \Pi_\omega(t,t)=I &t\in \R,
\endcases$$
where, for every $u_s\in H^1(\R^N)$, $\Pi_\omega(t,s)u_s$ is the
solution of \rff intro1.. with $u(s)=u_s$. \par We are interested in
the behaviour of the solutions of \rff intro1.. as $\omega
\rightarrow +\infty$. It is a well known fact that, given a Banach space ${\Cal M}$,  if a
function
$\sigma:\R \rightarrow {\Cal M}$ is almost periodic, the mean
value $$ \lim_{T\rightarrow +\infty}{{1}\over{2T}}\int_{-T}^T \sigma(p)\,dp =:\bar \sigma$$
exists. We observe that, for fixed $T>0$, $$\lim_{\omega \to
+\infty} \int_{-T}^T (\sigma(\omega p)-\bar \sigma)\, dp
=2T \lim_{\omega \to +\infty} {{1}\over{2\omega T}}\int_{-\omega
T}^{\omega T} (\sigma(p)-\bar \sigma)\, dp=0.$$  
Even if this
convergence is very weak, it suggests that the averaged equation
$$ u_t=\sum_{i,j=1}^N \bar a_{ij}\partial_i\partial_j u+ \bar
f(u)+ \bar g(x),\quad (t, x)\in\R_+\times\R^N\tag\dff intro2.. $$
should behave like a limit equation for \rff intro1.. as
$\omega \to +\infty$.\par Results of this kind have been known for quite a long time for
ordinary differential equations with almost periodic coefficients. For partial
differential equations, local results in this direction have been
obtained in an abstract setting (fit also for the study of
functional equations) by Hale and Verduyn Lunel \cite{\rfa
hallun..}. They consider an abstract semilinear parabolic equation
$$
u_t=Lu + f(\omega
p,u)+g(\omega p)\tag\dff intro3.. 
$$
in a Banach space $E$, where
$L$ is the generator of a strongly continuous semigroup of linear operators 
and $f(\cdot,u)$ and
$g(\cdot)$ are almost periodic.
They show the convergence of local solutions of \rff intro3.. to
solutions of the averaged equation 
$$
u_t=Lu+\bar f(u)+\bar g; \tag\dff intro4..
$$
moreover, they prove a continuation principle for strongly hyperbolic equilibria 
of \rff intro4.. and obtain an upper-semicontinuity result for local attractors
of the Poincar\'e map of \rff intro3... 

\par In a recent paper (\cite{\rff il..}), Ilyin
proposes a {\it global} criterion for comparison between the
process generated by \rff intro3.. and the semigroup generated by
the averaged problem \rff intro4... For autonomous equations like
\rff intro4.., it is well known that if $\bar f$ is dissipative and compact, then the
semiflow generated by \rff intro4.. possesses a compact global attractor in
$E$. In this case, it is possible to express the concept
of closeness of two semiflows in terms of the Hausdorff distance
of their attractors. As Ilyin shows in \cite{\rfa il..}, the same
can be done in the non-autonomous case. Ilyin considers an
abstract semilinear parabolic equation like \rff intro3.., where now
$L$ is a sectorial linear operator, and the corresponding averaged equation \rff intro4...
Using a notion of global attractor for families of processes 
introduced by Chepyzhov and
Vishik in \cite{\rfa chepvish..} (see Section 3), he shows that,
under suitable dissipativeness and compactness hypotheses, the
global attractor ${\Cal A}_\omega$ of \rff intro3.. converges in
the Hausdorff metric to the global attractor ${\Cal A}$ of \rff
intro4... Then he applies the abstract results to
reaction-diffusion, Navier-Stokes and damped wave equations on a
bounded domain $\Omega$.

The aim of our paper is to extend the results of \cite{\rfa il..} to
reaction-diffusion equations on the whole $\R^N$ with time
dependent principal part, like \rff intro1...
To this end, we cannot apply directly the abstract results of \cite{\rfa il..}. 
Indeed, since we are
working on the whole $\R^N$, the imbedding of
$H^1$ into $L^2$ is not compact; this makes much more difficult to
recover the asymptotic compactness of the processes generated by
\rff intro1... Even in the autonomous case, establishing the
existence of compact global attractors becomes then itself an
interesting task. In \cite{\rfa babinviscik2..} Babin and Vishik
overcame the difficulties arising from the lack of compactness by
introducing weighted Sobolev spaces. The choice of weighted
spaces, however, imposes some severe conditions on the forcing
term $g$ and on the initial data. Very recently, Wang (\cite{\rfa
wang..}) established the asymptotic $L^2$-compactness of the
semiflows and consequently the existence of global $(L^2-L^2)$
attractors for reaction-diffusion equations on $\R^N$ (or, more
generally, on unbounded subdomains of $\R^N$) avoiding the use of
weighted spaces. Following Wang's pattern, we shall prove uniform
asymptotic $L^2$-compactness of the processes generated by \rff
intro1... Then we shall obtain the asymptotic $H^1$-compactness by
a continuity argument similar to that of \cite{\rfa antopriz..}
and \cite{\rfa priz..}.\par On the other hand, since we assume that the
principal part is time-dependent, in the variation of constant
formula the linear semigroup $e^{-Lt}$ has to be replaced by the
linear processes $V_{\omega}(t,s)$ generated by the linear
equations $$u_t=\sum_{i,j=1}^N a_{ij}(\omega
t)\partial_i\partial_j u. \tag\dff intro6..$$ As a consequence, we have to prove
also the convergence of $V_{\omega}(t,s)$ to $e^{-\bar A t}$ as
$\omega \rightarrow +\infty$, where $e^{-\bar A t}$ denotes the
linear semigroup generated by the averaged linear equation. This
is done by mean of an explicit representation of the solutions of
the linear equations in terms of their Fourier transforms.

The paper is organized as follows. In Section 2 we introduce
notations and some necessary preliminaries; moreover, we obtain
some a-priori estimates for equation \rff intro1.. and we deduce
the existence of uniformly absorbing sets for the corresponding
process. In Section 3, we recall some basic properties of almost
periodic functions and the notion of uniform attractor
for a family of processes introduced by Chepyzhov and Vishik in \cite{\rfa chepvish..}; 
then we prove the existence of compact global uniform attractors for the families of
processes associated to \rff intro1... In Section 4 we investigate the behaviour of the
solutions of \rff intro1.. as $\omega \to +\infty$, proving that
the solutions of \rff intro1.. with initial datum $u_0\in
H^1(\R^N)$ converge, as $\omega \to +\infty$, to the solution of
\rff intro2.. with the same initial datum. Finally, we prove the upper-semicontinuity of the
family of the uniform attractors of \rff intro1.. as $\omega \to +\infty$, showing that the
uniform attractor of \rff intro1.. is $H^1$-close to that of \rff intro2..
for sufficiently large $\omega$. \par We would like to remark that the same
results hold for a family of reaction-diffusion equations of the
form $$u_t=\nu(\omega t)\sum_{i,j=1}^N \partial_i (a_{ij}(x)
\partial_j u)+f(\omega t,u)+g(\omega t ,x),\quad (t,x)\in
\R^+\times \Omega, $$ with Dirichlet or Neumann boundary
conditions on a bounded domain $\Omega\subset R^N$. To this end, it
suffices to replace the Fourier transform representation of the linear processes with their
spectral representation on a basis of eigenfunctions of the linear operator $\sum_{i,j=1}^N
\partial_i (a_{ij}(x)\partial_j )$ with the given boundary conditions.

\newsection\head \the\secnumber. Preliminaries\endhead

We consider the equation
$$
u_t=\sum_{i,j=1}^N a_{ij}(\omega t)\partial_i\partial_j u-a_0(\omega t)u+f(\omega t,u)+
g(\omega t,x),\quad x\in\R^N,\tag\dff equazione..
$$
where $\omega$ is a positive constant. \par We make the following assumptions: 
the functions $a_{ij}$ and $a_0$ are H\"older continuous on $\R$ with exponent $\theta$,
$a_{ij}(\tau)=a_{ji}(\tau)$ for $i,j=1$, \dots, $N$ and for all $\tau\in\R$,
and there exist positive constants $\nu_1\geq\nu_0>0$ and $C>0$ such that
$$
\nu_0|\xi|^2\leq\sum_{i,j=1}^N a_{ij}(\tau)\xi_i\xi_j\leq\nu_1|\xi|^2\quad\text{for all
$\tau\in\R$ and
$\xi\in\R^N$}\tag\dff a1..
$$
and
$$
|a_0(\tau)|\leq C\quad\text{for all $\tau\in\R$.}
$$
Moreover,
$$
\|g(\tau,\cdot)\|_{ L^2}\leq C\quad\text{for all $\tau\in\R$}\tag\dff a2..
$$
and there exist $g_0\in L^2(\R^N)$ and $0<\theta\leq 1$ such that
$$
|g(\tau_1,x)-g(\tau_2,x)|\leq g_0(x)|\tau_1-\tau_2|^\theta\quad\text{for all $\tau_1,\tau_2\in\R$
and for a.e. $x\in\R^N$.}\tag\dff a3..
$$
Finally,
$$
f(\tau,0)= 0,\quad|f_u(\tau,u)|\leq C(1+|u|^\beta)\quad\text{for all
$u,\tau\in\R$}\tag\dff a4..
$$
and
$$
|f(\tau_1,u)-f(\tau_2,u)|\leq C(|u|+|u|^{\beta+1})|\tau_1-\tau_2|^\theta
\quad\text{for all $u,\tau_1,\tau_2\in\R$,}\tag\dff a5..
$$
where
$$
0\leq\beta\quad\text{if $N\leq2$};\quad 0\leq\beta\leq 2^*/2-1\quad\text{if $N\geq 3$.}
\tag\dff a6..
$$

For $t\in\R$ and $\omega>0$ we define the operator $A_\omega(t)\colon H^2(\R^N)\to L^2(\R^N)$
by
$$
A_\omega(t)u:=-\sum_{i,j=1}^N a_{ij}(\omega t)\partial_i\partial_j u,\quad u\in H^2(\R^N).
$$
Then $A_\omega(t)$ is a self-adjoint positive operator in $L^2(\R^N)$ and our assumptions on the
coefficients $a_{ij}(\tau)$ imply that the abstract parabolic
equation
$$
\dot u=-A_\omega(t)u
$$
generates a linear process
$$
U_\omega(t,s)\colon L^2(\R^N)\to L^2(\R^N), \quad t\geq s,
$$
such that
$$
\|U_\omega(t,s)u\|_{L^2}\leq M\|u\|_{L^2},\quad u\in L^2(\R^N),\tag\dff expest1..
$$
$$
\|U_\omega(t,s)u\|_{H^1}\leq M\|u\|_{H^1},\quad u\in H^1(\R^N),\tag\dff expest2..
$$
and
$$
\|U_\omega(t,s)u\|_{H^1}\leq M(1+(t-s)^{-1/2})\|u\|_{L^2},\quad u\in L^2(\R^N),
\tag\dff expest3..
$$
where $M$ is a positive constant
(see e.g. \cite{\rfa Pazy.., Ch.5}, \cite{\rfa tana..}). \par A useful
explicit representation of
$U_\omega(t,s)$ can be given in terms of its Fourier transform. We denote by ${\Cal
F}v$ the Fourier-Plancherel transform of $v\in L^2(\R^N)$, normalized in such a way that, for
$v\in L^2(\R^N)\cap L^1(\R^N)$, 
$$
({\Cal F}v)(\xi)={1\over{(2\pi)^{N/2}}}\int_{\R^N}e^{-i x\cdot\xi}v(x)\, dx.
$$
It is well known that ${\Cal F}$ is an isometry of $L^2(\R^N)$ onto itself, and $v\in
H^k(\R^N)$ if and only if
$(1+|\xi|^2)^{k/2}({\Cal F}v)(\xi)\in L^2(\R^N)$. Moreover,
$$
\|u\|_{H^1}^2=\int_{\R^N}(1+|\xi|^2)({\Cal F}u)(\xi)^2\, d\xi.\tag\dff norm..
$$
Then an easy computation gives
$$
({\Cal F}(U_\omega(t,s)u))(\xi)=\exp\left\{-\int_s^t(\sum_{i,j=1}^N a_{ij}(\omega
p)\xi_i\xi_j)\,dp\right\}({\Cal F}u)(\xi).\tag\dff rappr..
$$
An immediate consequence of \rff rappr.. is that the constant $M$ in \rff expest1.. -- \rff
expest3.. depends only on $\nu_0$, so in particular is independent of $\omega$.

As for the nonlinear term, conditions \rff a4.. and \rff a5.. and the Sobolev embedding
Theorem imply that the Nemitski operator
$$
\hat f\colon \R\times H^1(\R^N)\to L^2(\R^N)
$$
is well defined and satisfies
$$
\|\hat f(\tau, u)\|_{L^2}\leq \tilde C(\|u\|_{L^2}+\|u\|_{H^1}^{\beta+1}),\quad\tau\in\R, u\in
H^1(\R^N)\tag\dff bound..
$$
and
$$
\multline
\|\hat f(\tau_1,u_1)-\hat f(\tau_2,u_2)\|_{L^2}
\leq \tilde C(1+\|u_1\|_{H^1}^{\beta+1}+\|u_2\|_{H^1}^{\beta+1})|\tau_1-\tau_2|^\theta\\
+\tilde C(1+\|u_1\|_{H^1}^{\beta}+\|u_2\|_{H^1}^{\beta})\|u_1-u_2\|_{H^1},\quad
\tau_1,\tau_2\in\R,\quad u_1,u_2\in H^1(\R^N),
\endmultline\tag\dff lip..
$$
where $\tilde C$ is a positive constant depending only on $C$, $\nu_0$, $\nu_1$ and $\beta$.
By classical results of \cite{\rfa fried..}, \cite{\rfa He..} and \cite{\rfa Pazy..}, for every
$s\in\R$ and for every $u_s\in H^1(\R^N)$ the semilinear Cauchy
problem
$$
\cases
\dot u=-A_\omega(t)u-a_0(\omega t)u+\hat f(\omega t,u)+g(\omega t)\\
u(s)=u_s
\endcases
\tag\dff pbcauchy..
$$
is locally well-posed and hence possesses a unique maximal classical solution $u\in
C^0([s,s+T[, H^1)\cap C^1(]s,s+T[,L^2)$, $T$ depending on $s$ and $u_s$. Moreover,  $u$
satisfies the variation of constant formula
$$
u(t)=U_\omega(t,s)u_s+\int_s^t U_\omega(t,p)
(-a_0(\omega p)u(p)+\hat f(\omega p,u(p))+g(\omega p))\,dp,\quad t\geq s.
$$

The following set of dissipativeness and monotonicity conditions ensures that the solutions
of \rff pbcauchy.. are global and bounded:
$$
a_0(\tau)\geq\lambda_0>0\quad\text{for all $\tau\in\R$};\tag\dff a7..
$$
$$
 f(\tau,u)u\leq 0;\quad f_u(\tau,u)\leq L\quad\text{for all $u,\tau\in\R$}.
\tag\dff a8..
$$
We start with the following a-priori estimates in  $L^2$:
\proclaim{Lemma \dft priori1..}
Let $u_\omega \co [s,s+T[ \to H^1(\R^N)$ be the maximal solution of the Cauchy problem \rff
pbcauchy... If $\| u_s\|_{L^2}\leq R$, then, for $t\in[s,s+T[$,
$$
\|u_\omega(t)\|_{L^2}^2\leq e^{-\lambda_0 (t-s)}R^2+{{C^2}\over{\lambda_0^2}}
$$
\endproclaim
\demo{Proof}
For $t\in ]s,s+T[$, we have
$$
\alignedat1
&{{d}\over{dt}}{1\over2}\|u_\omega(t)\|_{L^2}^2
=\langle u_\omega(t),\dot u_\omega(t)\rangle\\
&=\langle u_\omega(t),-A_\omega(t) u_\omega(t)-a_0(\omega t) u_\omega(t)+
\hat f(\omega t,u_\omega(t))+g(\omega t)\rangle\\
&\leq-\langle u_\omega(t),A_\omega(t)u_\omega(t)\rangle-
\lambda_0\langle u_\omega(t),u_\omega(t)\rangle+
\langle u_\omega(t),\hat f(\omega t,u_\omega(t))\rangle+\langle u_\omega(t),g(\omega t)\rangle.
\endalignedat
$$
By \rff a8.. and by Young's inequality, we get
$$
\multline
{{d}\over{dt}}{1\over2}\|u_\omega(t)\|_{L^2}^2+\langle u_\omega(t),A_\omega(t)u_\omega(t)\rangle+
\lambda_0\|u_\omega(t)\|_{L^2}^2\\
\leq \langle u_\omega(t),g(\omega t)\rangle\leq
{\lambda_0\over 2}
\|u_\omega(t)\|_{L^2}^2+{1\over{2\lambda_0}}\|g(\omega t)\|_{L^2}^2.
\endmultline
$$
By \rff a2.. it follows that
$$
{{d}\over{dt}}\|u_\omega(t)\|_{L^2}^2+\lambda_0\|u_\omega(t)\|_{L^2}^2\leq
{{C^2}\over
\lambda_0}.
$$
Multiplication by $e^{\lambda_0 t}$ and integration yields
$$
\|u_\omega(t)\|_{L^2}^2\leq e^{-\lambda_0 (t-s)}\|u_\omega(s)\|_{L^2}^2+{{C^2}\over
{\lambda_0^2}},\tag\dff prio1..
$$
and the conclusion follows.
\qed\enddemo

In order to get $H^1$-estimates, we need the following lemmas:

\proclaim{Lemma \dft appr..}
Let $u\in H^2(\R^N)$. Then $\langle\hat f(\omega t,u),-\Delta u\rangle\leq
L \langle u,-\Delta u\rangle$ for all $t\in\R$, where $L$ is the constant of condition \rff
a8...
\endproclaim
\demo{Proof}
For $n\in\N$, choose a function $h_n\in C^\infty(\R)$, with $0\leq h'_n(u)
\leq 1$ for all $u\in\R$, such that
$$
h_n(u)=\cases
u&\text{if $-n\leq u\leq n$}\\
n+1&\text{if $2n\leq u$}\\
-(n+1) &\text{if $ u\leq -2n$}
\endcases
$$
Let us fix $t\in\R$ and define $f_n(\omega t,u):=f(\omega t,h_n(u))$. By \rff a8.., it follows
that
$f_n(\omega t,0)=0$,
$|(f_n)_u(\omega t,u)|$ is bounded on $\R$ and $(f_n)_u(\omega t,u)\leq L$ for all $u\in\R$.
By Proposition IX.5 in
\cite{\rfa bre..}, it follows that $f_n(\omega t,u(\cdot))\in H^1(\R^N)$ and
$\nabla f_n(\omega t,u(\cdot))
=(f_n)_u(\omega t,u(\cdot))\cdot\nabla u$. Then, for all $n\in\N$,  we have
$$
\multline
\langle\hat f_n(\omega t,u),-\Delta u\rangle
=\int_{\R^N} (f_n)_u(\omega t, u(x))|\nabla u(x)|^2
\, dx\\
\leq L\int_{\R^N} |\nabla u(x)|^2
\, dx=L\langle u,-\Delta u\rangle
\endmultline
$$
The proof will be complete if we show that $f_n(\omega t,u(\cdot))\to
f(\omega t,u(\cdot))$ in $L^2(\om)$
as $n\to\infty$.
Now, since $f_n(\omega t, u(x))\to f(\omega t,u(x))$ almost everywhere in $\R^N$  and
the estimates
$$
|f_n(\omega t,u(x))|\leq C(|u(x)|+|u(x)|^{\beta+1}),
$$
$$
|f(\omega t,u(x))|\leq C(|u(x)|+|u(x)|^{\beta+1})
$$
hold, the conclusion follows from the Lebesgue dominated convergence theorem.
\qed\enddemo

\proclaim{Lemma \dft fourier..}
For all $u\in H^2(\R^N)$ and for all $t\in\R$
$$
-\langle A_\omega(t)u,\Delta u\rangle\geq \nu_0 \|\Delta u\|_{L^2}^2.
$$
\endproclaim
\demo{Proof}
Again denoting by ${\Cal F}v$ the
Fourier-Plancherel transform of $v\in L^2(\R^N)$, we have
$$
\multline
-\langle A_\omega(t)u,\Delta u\rangle
=\langle{\Cal F}( A_\omega(t)u),{\Cal F}(-\Delta u)\rangle\\
=\int_{\R^N}(\sum_{i,j=1}^N a_{ij}(\omega t)\xi_i\xi_j)({\Cal
F}u)(\xi)(\sum_{l=1}^N\xi_l^2)({\Cal F}u)(\xi)\,d\xi\\
\geq\nu_0\int_{\R^N}\left(\sum_{l=1}^N\xi_l^2 ({\Cal F}u)(\xi)\right)^2\,d\xi
=\nu_0\|\Delta u\|^2_{L^2}.
\endmultline
$$
\qed\enddemo

Now we are able to prove

\proclaim{Lemma \dft priori2..}
Let $u_\omega \co [s,s+T[ \to H^1(\R^N)$ be the maximal solution of the Cauchy problem \rff
pbcauchy...
There exist two positive constants $K_1$ and $K_2$, depending only on
$C$, $\nu_0$, $\nu_1$, $\lambda_0$ and $L$, such that,
if $\|u_s\|_{H^1}\leq R$, then, for $t\in[s,s+T[$,
$$
\|u_\omega(t)\|_{H^1}^2\leq K_1R^2e^{-\lambda_0 (t-s)}+K_2
$$
\endproclaim
\demo{Proof} For $t\in ]s,s+T[$, by Lemma \rft fourier.. and by \rff a7..  we have
$$
\alignedat1
&{{d}\over{dt}}{1\over2}\|\nabla u_\omega(t)\|^2_{L^2}=
-\langle \Delta u_\omega(t),\dot u_\omega(t)\rangle\\
&=\langle-\Delta  u_\omega(t), -A_\omega(t) u_\omega (t)-a_0(\omega t) u_\omega(t) +\hat
f(\omega t,u_\omega(t))+g(\omega t)\rangle\\
& \leq-\nu_0\|\Delta u_\omega(t)\|_{L^2}^2-\lambda_0 \|\nabla u_\omega(t)\|_{L^2}^2-
\langle \Delta u_\omega(t),\hat f(\omega t,u_\omega(t))\rangle-
\langle \Delta u_\omega(t),g(\omega t)\rangle.
\endalignedat
$$
By \rff a2.., by Lemma \rft appr.. and by Young's inequality we obtain
$$
{{d}\over{dt}}\|\nabla u_\omega(t)\|^2_{L^2}\leq
-\nu_0\|\Delta u_\omega(t)\|_{L^2}^2-2(\lambda_0-L) \|\nabla u_\omega(t)\|_{L^2}^2+
{{C^2}\over{\nu_0}}.\tag\dff trg1..
$$
Let $\delta>0$ and let $v\in H^2(\R^N)$. We have
$$
\|\nabla v\|_{L^2}^2\leq {\delta\over2}\|\Delta v\|_{L^2}^2+{1\over{2\delta}}\|v\|_{L^2}^2,
$$
whence
$$
-\|\Delta v\|_{L^2}^2\leq -{2\over\delta}\|\nabla v\|_{L^2}^2+{1\over{\delta^2}}\|v\|_{L^2}^2.
\tag\dff trg2..
$$
By \rff trg1.. and \rff trg2.., choosing  $\delta:=\nu_0/L$, we obtain
$$
{{d}\over{dt}}\|\nabla u_\omega(t)\|^2_{L^2}\leq
-2\lambda_0 \|\nabla u_\omega(t)\|^2_{L^2} +{{L^2}\over{\nu_0}}
\|u_\omega(t)\|_{L^2}^2+{{C^2}\over{\nu_0}}.
$$
By \rff prio1..,
$$
{{d}\over{dt}}\|\nabla u_\omega(t)\|^2_{L^2}+2\lambda_0 \|\nabla u_\omega(t)\|^2_{L^2}
\leq{{L^2}\over{\nu_0}}\|u_\omega(s)\|_{L^2}^2
e^{-\lambda_0 (t-s)}
+\left({{L^2}\over{\nu_0\lambda_0^2}}+{1\over{\nu_0}}\right)C^2.
$$
Multiplication by $e^{2\lambda_0 t}$ and integration yields
$$
\multline
\|\nabla u_\omega(t)\|_{L^2}^2\\
\leq e^{-2\lambda_0(t-s)}\|\nabla u_\omega(s)\|_{L^2}^2
+{{L^2}\over{\nu_0\lambda_0}}\|u_\omega(s)\|_{L^2}^2e^{-\lambda_0(t-s)}+
\left({{L^2}\over{2\nu_0\lambda_0^3}}+{1\over{\lambda_0\nu_0}}\right)C^2,
\endmultline
$$
and the conclusion follows.
\qed\enddemo
As a consequence, we have the following result:
\proclaim{Proposition \dft absorb1..}
Let $u_\omega \co [s,s+T[ \to H^1(\R^N)$ be the maximal solution of the Cauchy problem \rff
pbcauchy... Then
\roster
\item $T=+\infty$;
\item if $\|u_s\|_{H^1} \leq R$, then, for every $t\geq s$,
$\|u_\omega(t)\|_{H^1}^2\leq K_1R^2+K_2$, where $K_1$ and $K_2$
are independent of $R$ and $\omega$;
\item there exists a positive constant $K$, and for every $R>0$ there exists
$T(R)>0$ such that, whenever $\|\cl u_s\|_{H^1} \leq R$,
$\|u_\omega(t)\|_{H^1}< K$ for all $t$ such that $(t-s)\geq T(R)$. Both $K$ and $T(R)$, besides
$R$,  depend only on $C$, $\nu_0$, $\nu_1$, $\lambda_0$
and $L$. In particular, they are independent of $s$ and $\omega$.
\endroster \qed
\endproclaim

Proposition \rft absorb1.. says also that the global process generated by \rff pbcauchy..
possesses a bounded absorbing set in $H^1(\R^N)$ independent of $\omega$. \par  
We end this section with a result which will be useful in proving the asymptotic compactness
of the processes generated by \rff pbcauchy... 

\proclaim{Lemma \dft ball1..}
Let $u_\omega \co [s,+\infty[ \to H^1(\R^N)$ be the solution of the Cauchy problem \rff
pbcauchy.., with $\|u_s\|_{H^1} \leq R$. Assume moreover that the set
$\{\,g(\tau,\cdot)\mid\tau\in\R\,\}$ is compact in $L^2(\R^N)$.
Then, for every
$\eta>0$, there exist two positive constants $\cl k(R)$ and $\cl T(R)$ such that, if
$(t-s)\geq\cl T(R)$ and $k
\geq
\cl k(R)$,
$$
\int_{\{|x|> k\}}|u_\omega(t,x)|^2 \, d x \leq \eta.
$$
The constants $\cl k(R)$ and $\cl T(R)$, besides $R$ and $\eta$, depend only
on $C$, $\nu_0$, $\nu_1$, $\lambda_0$
and $L$.
In particular, they
are independent of $s$ and $\omega$.
\endproclaim
\demo{Proof}
We adapt to the non-autonomous case the proof of Lemma 5 in \cite{\rfa wang..}, being careful
that all the estimates involved are independent of $\omega$.\par  Let $\theta\co \R_+ \to
\R$ be a smooth function such that $0\leq \theta(s) \leq 1$ for $s
\in \R_+$, $\theta(s)=0$ for $0\leq s \leq 1$ and $\theta(s)=1$ for $s \geq 2$.
Let $D := \sup_{s \in \R_+} |\theta'(s)|$. For $k \in
\N$, let us define the multiplication operator
$$\Theta_k\co H^1(\R^N) \to H^1(\R^N), \quad
(\Theta_ku)(x):=\theta(|x|^2/k^2)u(x).
$$
By \rff a7.. and \rff a8.., we have
$$
\alignedat1
&{{d}\over{dt}}{1\over2}\int_{\R^N}
\theta(|x|^2/k^2)|u_\omega(t,x)|^2\, d x =
{{d}\over{dt}}{1\over2}\langle
\Theta_k u_\omega(t), u_\omega(t)\rangle
=\langle\Theta_k u_\omega(t),\dot u_\omega(t)\rangle\\
& =\langle \Theta_k u_\omega(t),-A_\omega(t)
u_\omega(t)-a_0(\omega t) u_\omega(t)+\hat f(\omega t,u_\omega(t))+g(\omega t)\rangle\\
&\leq\langle \Theta_k u_\omega(t),-A_\omega(t)u_\omega(t)\rangle-\lambda_0\langle \Theta_k
u_\omega(t),u_\omega(t)\rangle+\langle \Theta_k u_\omega(t),g(\omega t)\rangle
\endalignedat
$$
and hence
$$\multline
{{d}\over{dt}}\langle \Theta_k u_\omega(t), u_\omega(t)\rangle +2 \lambda_0\langle \Theta_k
u_\omega(t),u_\omega(t)\rangle\\
 \leq - 2\langle \Theta_k
u_\omega(t),A_\omega(t)u_\omega(t)\rangle +2 \langle \Theta_k u_\omega(t),g(\omega t)\rangle.
\endmultline
$$
Since
$$
\alignedat1 &\langle \Theta_k
u_\omega(t),A_\omega(t)u_\omega(t)\rangle= \int_{\R^N}
\theta(|x|^2/k^2) \sum_{i,j=1}^N a_{ij}(\omega t)\partial_i
u_\omega(t,x)\partial_ju_\omega(t,x)\, d x
\\ &+ \int_{\R^N} \theta'(|x|^2/k^2)
u_\omega(t,x) {2 \over {k^2}}\sum_{i,j=1}^N a_{ij}(\omega t) x_i \partial_j u_\omega(t,x)\, d x,
\endalignedat
$$
it follows that
$$\alignedat1
&-\langle \Theta_k
u_\omega(t),A_\omega(t)u_\omega(t)\rangle \leq \int_{\R^N} \theta'
(|x|^2/k^2) |u_\omega(t,x)|| \nabla u_\omega(t,x)|{2 \over {k^2}}|x|\, d x \\
&\leq 2 D \int_{ \{k \leq |x| \leq \sqrt2 k\}}{{|x|} \over{k^2}}|u_\omega(t,x)|
|\nabla u_\omega(t,x)| \, d x  \\
&\leq {{2 \sqrt2 D}\over k} \int_{\{k \leq |x| \leq \sqrt2 k\}}
|u_\omega(t,x)| |\nabla u_\omega(t,x)| \, d x  \\
&\leq {{2 \sqrt2 D}\over k} \|u_\omega(t)\|_{L^2}\|\nabla u_\omega(t) \|_{L^2}.
\endalignedat
$$
So, by Proposition \rft absorb1.., for $ (t-s) \geq T(R)$, we have
$$
- \langle \Theta_k
u_\omega(t),A_\omega(t)u_\omega(t)\rangle\leq {{2 \sqrt2 D K^2}\over k}
$$
Let $\eta >0$ and choose $k=k(\eta)$ such that
$${{2 \sqrt2 D K^2}\over k}< \eta.$$
Then for $(t-s) > T(R)$ and $k> k(\eta)$,
we obtain
$$
{{d}\over{dt}}\langle \Theta_k u_\omega(t), u_\omega(t)\rangle +2 \lambda_0\langle \Theta_k
u_\omega(t),u_\omega(t)\rangle
\leq 2 \eta + 2 \langle \Theta_k u_\omega(t),g(\omega t)\rangle.
$$
By Young's inequality, we have
$$
\langle \Theta_k u_\omega(t),g(\omega t)\rangle \leq {{\lambda_0}\over 2} \langle \Theta_k
u_\omega(t),u_\omega(t) \rangle + {1 \over {2 \lambda_0}} \int_{\R^N}
\theta(|x|^2/k^2)  g(\omega t,x)^2 \,d x .
$$
Since we have assumed that $\{\,g(\tau,\cdot)\mid\tau\in\R\,\}$ is compact in $L^2(\R^N)$,
there exists
$k'=k'(\eta)$ such that, if $k>k'(\eta)$,
$$
{1 \over {2 \lambda_0}} \int_{\R^N}
\theta(|x|^2/k^2)  g(\omega t,x)^2 \, d x  \leq \eta\quad\text{for all $t\in\R$;}
$$
As a consequence, for $(t-s)>T(R)$ and for $k>\max \{k(\eta), k'(\eta)\}$
$$ {{d}\over{dt}}\langle \Theta_k u_\omega(t), u_\omega(t)\rangle + \lambda_0\langle \Theta_k
u_\omega(t),u_\omega(t)\rangle
\leq 4 \eta.
$$
Multiplication by $e^{\lambda_0 t}$ and integration yields
$$
e^{\lambda_0 t} \langle \Theta_k u_\omega(t), u_\omega(t)\rangle - e^{\lambda_0(s+ T(R))}
\langle
\Theta_k u_\omega(s+T(R)), u_\omega(s+T(R))\rangle
\leq {{4\eta}\over{\lambda_0}} e^{\lambda_0 t}
$$
for $(t-s) >T(R)$.
It follows that, for $(t-s)>T(R)$,
$$
\multline
\langle \Theta_k u_\omega(t), u_\omega(t)\rangle
\leq e^{- \lambda_0 ((t-s)-T(R))}\langle \Theta_k
u_\omega(s+T(R)), u_\omega(s+T(R))\rangle + {{4\eta}\over\lambda_0} \\
\leq e^{- \lambda_0 ((t-s)-T(R))}
K^2 + {{4\eta}\over\lambda_0}.
\endmultline
$$
Finally, for $(t-s)\geq T(R) + \lambda_0^{-1}\log(\eta^{-1})$ and for
$k>\max\{k(\eta),k'(\eta)\}$, we get
$$
\int_{\{|x|>\sqrt2 k\}}|u_\omega(t,x)|^2 \, d x \leq \int_{\R^N}
\theta(|x|^2/k^2)| u_\omega(t,x)|^2 \, d x \leq \left( K^2 + {4
\over \lambda_0} \right) \eta,
$$
and the proof is complete.
\qed\enddemo

\newsection\head\the\secnumber. Existence of the compact global attractors\endhead

It is well known (see e.g. \cite{\rfa Lady..}, \cite{\rfa hala..}, \cite{\rfa babinviscik..}
and \cite{\rfa temam..}) that if a continuous semigroup
$P(t)$, acting on a complete connected  metric space $X$, is bounded, pointwise-dissipative
and asymptotically compact, then it possesses a compact global attractor. The attractor is
non-empty, connected, strictly invariant, and can be characterized as the union of all
complete bounded trajectories of $P(t)$. This idea can be quite naturally extended to the
class of processes generated by periodically time-dependent partial differential equations,
since such systems undergo a discrete semigroup structure given by the period map. In more
general situations, like the almost periodic case considered here, the leading property of
invariance fails and new approaches had to be developed.\par
In \cite{\rfa hara..}, Haraux proposed a notion of attractor for a process $\Pi(t,s)$ based
on the concept of minimality rather than invariance. However, as it was suggested by the same
Haraux, the theory of skew-product flows, at the expense of introducing an extended phase
space, provides the right extension of invariance. This alternative approach, developed by
Chepyzhov and Vishik in \cite{\rfa chepvish..}, turns out to be particularly well suited
if the process is generated by an almost periodic partial differential equation. \par
We shall describe this approach in the context of equation \rff equazione...\par  

We define ${\Cal M}_1$ as the space of $N\times N$ real symmetric matrices and ${\Cal
M}_2:=\R$; moreover, we denote by ${\Cal M}_3$ the set
$$
{\Cal M}_3:=\{\,\Psi\colon\R\to\R\mid \Psi(0)=0, \,\|\Psi\|_{{\Cal M}_3}<+\infty\,\},
$$
where
$$
\|\Psi\|_{{\Cal M}_3}:=\sup_{u\in\R}{{|\Psi_u(u)|}\over{1+|u|^{\beta}}}.
$$
Finally, we set ${\Cal M}_4:=L^2(\R^N)$. \par \medskip

Besides conditions \rff a1..--\rff a6.. and \rff a7..--\rff a8.., from now on we assume that
also the following condition is satisfied:

\smallskip

\item\item{\bf (AP)} the functions $t\mapsto (a_{ij}(t))_{ij}\in {\Cal M}_1$;
$t\mapsto a_0(t)\in{\Cal M}_2$; $t\mapsto f(t,\cdot)\in {\Cal M}_3$ and $t\mapsto g(t,\cdot)
\in{\Cal M}_4$ are almost periodic.

\smallskip

By Bochner's criterion (see e.g.
\cite{\rfa levzhik..}), whenever
$\sigma\colon \R\to{\Cal M}$ is almost periodic, the set of all translations
$\{\,\sigma(\cdot+h)\mid h\in\R\,\}$ is precompact in
$C_b(\R,{\Cal M})$. The closure of this set in $C_b(\R,{\Cal M})$ is called the hull of $\sigma$
and is usually denoted by ${\Cal H}(\sigma)$; if $\zeta\in {\Cal H}(\sigma)$, then
$\zeta$ is almost periodic and ${\Cal H}(\zeta)={\Cal H}(\sigma)$. \par 
For an almost periodic function $\sigma$, the mean value
$$
\lim_{T\to+\infty}{1\over{2T}}\int_{-T}^{T}\sigma(t)\, dt=:\bar\sigma\in{\Cal M}
$$
exists. More remarkably, (see again \cite{\rfa levzhik..}) there exists a
bounded decreasing function
$\mu\colon \R_+\to\R_+$, $\mu(T)\to 0$ as $T\to\infty$, such that
$$
\|(1/T)\int_s^{s+T}(\zeta(t)-\bar\sigma)\,dt\|_{{\Cal M}}\leq\mu(T)\quad
\text{for all $s\in\R$ and all $\zeta \in {\Cal H}(\sigma)$.}\tag\dff media..
$$

We will denote by $\Sigma_1$, $\Sigma_2$, $\Sigma_3$ and $\Sigma_4$ the hulls
of the functions $t\mapsto (a_{ij}(t))_{ij}$,
$t\mapsto a_0(t)$, $t\mapsto f(t,\cdot)$ and $t\mapsto g(t,\cdot)$ in $C_b(\R,{\Cal M}_1)$,
$C_b(\R,{\Cal M}_2)$, $C_b(\R,{\Cal M}_3)$ and $C_b(\R,{\Cal M}_4)$ respectively.
The corresponding mean values will be denoted by $(\bar a_{ij})\in{\Cal M}_1$, $\bar
a_0\in{\Cal M}_2$, $\bar f(\cdot)\in {\Cal M}_3$ and $\bar g(\cdot)\in{\Cal M}_4$. Moreover,
let us set $\Sigma:=\Sigma_1\times\Sigma_2\times\Sigma_3\times\Sigma_4$.

\remark{Remark} It is an easy exercise
to check that properties \rff a1..--\rff a6.. and \rff a7..--\rff a8.. are satisfied
by any element of $\Sigma_1$, $\Sigma_2$, $\Sigma_3$ and $\Sigma_4$, as well as by
the corresponding mean values (with the same
constants!). Hence, the results of Lemmas \rft priori1.. and
\rft priori2.., of
Proposition \rft absorb1.. and of Lemma \rft ball1.. still hold true if we replace
$a_{ij}(\tau)$,
$a_0(\tau)$, $f(\tau,u)$ and $g(\tau,x)$ with arbitrary elements of $\Sigma_1$, $\Sigma_2$,
$\Sigma_3$ and $\Sigma_4$ or with the corresponding mean values.
\endremark

\medskip

\remark{Remark}
Applying \rff media.. to $f(\tau,u)$, we have
$$
|(1/T)\int_s^{{s+T}}(f(\tau,u)-\bar f(u))\, d\tau|\leq \mu(T)(|u|+|u|^{\beta+1})\quad
\text{for all $s$ and $u\in\R$,}
$$
and integration yields
$$
\|(1/T)\int_s^{{s+T}}(\hat f(\tau,u)-\Hat{\Bar f}(u))\, d\tau\|_{L^2}\leq
K\mu(T)(\|u\|_{L^2}+\|u\|_{H^1}^{\beta+1}).\tag\dff 31..
$$
\endremark

\medskip

As a consequence of
Lemmas \rft priori1.. and
\rft priori2.. and of Proposition \rft absorb1.., for any
$\sigma=((\alpha_{ij}),\alpha_0,\phi,\gamma)\in\Sigma$ and for any $\omega>0$, the equation
$$
u_t=\sum_{i,j=1}^N \alpha_{ij}(\omega t)\partial_i\partial_j u-\alpha_0(\omega t)u+
\phi(\omega t,u)+
\gamma(\omega t,x),\quad x\in\R^N
$$
generates a global process $\Pi^\sigma_\omega(t,s)$ in the space $H^1(\R^N)$. \par
According to \cite{\rfa chepvish..}, now we are able to give the
following \proclaim{Definition \dft attractor.. (Chepyzhov and
Vishik, '94)} A closed set ${\Cal A}^\Sigma_\omega$ is said to be
the $\Sigma$-uniform attractor of the family of processes
$\{\,\Pi^\sigma_\omega\mid\sigma\in\Sigma\,\}$ iff \roster
\item For every bounded set $B\subset H^1(\R^N)$
$$
\lim_{t\to\infty}\sup_{\sigma\in\Sigma}\dist(\Pi^\sigma_\omega(t,s)B,{\Cal A}^\Sigma_\omega)=0
\quad\text{for all $s\in\R$;}
$$
\item ${\Cal A}^\Sigma_\omega$ is minimal among all closed subsets of $H^1(\R^N)$ satisfying
property (1), i.e.
${\Cal A}^\Sigma_\omega\subset{\Cal A}'$ for every closed set
${\Cal A}'\subset H^1(\R^N)$ which satisfies property (1).
\endroster
\endproclaim

Our first goal is to prove that the almost periodic dissipative equation
\rff equazione.. possesses a $\Sigma$-uniform attractor in $H^1(\R^N)$. Following \cite{\rfa
chepvish..}, we introduce the extended phase-space $\Sigma\times H^1(\R^N)$; for
$\omega>0$, we define on $\Sigma$ the unitary group of translations
$$
(T_\omega(h)\sigma)(\cdot):=\sigma(\cdot+\omega h).
$$
One can easily prove the following translation identity:
$$
\Pi^\sigma_\omega(t+h,s+h)=\Pi^{T_\omega(h)\sigma}_\omega(t,s), \quad h\in\R. \tag\dff
transliden..
$$
Thanks to \rff transliden.., we can associate to the family of processes
$\{\,\Pi^\sigma_\omega\mid\sigma\in\Sigma\,\}$ a (nonlinear) semigroup
$P_\omega(t)$ acting on the extended phase-space
$\Sigma\times H^1(\R^N)$, by the formula
$$
P_\omega(t)(\sigma,u):=(T_\omega(t)\sigma,\Pi_\omega^\sigma(t,0)u).
$$
 In \cite{\rfa chepvish..}, Chepyzhov and Vishik proved that, if the semigroup
$P_\omega(t)(\sigma,u)$ above is continuous, bounded, pointwise-dissipative and asymptotically
compact (and hence possesses a compact global attractor ${\Cal M}^\Sigma_\omega$), then the
projection of ${\Cal M}^\Sigma_\omega$ onto $H^1(\R^N)$ is the global $\Sigma$-uniform attractor
of the family of processes $\{\,\Pi^\sigma_\omega\mid\sigma\in\Sigma\,\}$. \par
Let us describe in some detail the results of \cite{\rfa chepvish..}.

\proclaim{Definition \dft fullbounded..}
A curve $t\to u(t)\in H^1(\R^N)$, $t\in\R$ is said to be a full solution of the process
$\Pi^\sigma_\omega(t,s)$ iff
$$
\Pi^\sigma_\omega(t,s)u(s)=u(t)\quad\text{for all $t\geq s$, $s\in\R$.}
$$
\endproclaim

\proclaim{Definition \dft kernel..}
The kernel of the process $\Pi^\sigma_\omega(t,s)$ is by definition the set ${\Cal
K}^\sigma_\omega$ of all full bounded solutions of the process $\Pi^\sigma_\omega(t,s)$.
We call the set
$$
{\Cal K}^\sigma_\omega(s):=\{\,u(s)\mid u(\cdot)\in {\Cal K}^\sigma_\omega\,\}\subset H^1(\R^N)
$$
the kernel section at time $s$.
\endproclaim

We introduce also the two projectors
$J_1$ and $J_2$ from $\Sigma\times H^1(\R^N)$ onto $\Sigma$ and $H^1(\R^N)$ respectively:
$J_1(\sigma, u):=\sigma$, $J_2(\sigma,u):=u$. Then we have

\proclaim{Theorem \dft ChepVish..\, (Chepyzhov and Vishik, '94)}
Assume that the semigroup $P_\omega(t)$
is continuous, bounded, pointwise-dissipative and asymptotically
compact, so it possesses a compact global attractor ${\Cal M}^\Sigma_\omega$.
Then
\roster
\item $J_2{\Cal M}^\Sigma_\omega=:{\Cal A}^\Sigma_\omega$ is the global
$\Sigma$-uniform attractor of the family of processes
$\{\,\Pi^\sigma_\omega\mid\sigma\in\Sigma\,\}$;
\item $J_1{\Cal M}^\Sigma_\omega=\Sigma$;
\item ${\Cal M}^\Sigma_\omega=\bigcup_{\sigma\in\Sigma}\,\{\sigma\}\times {\Cal
K}^\sigma_\omega(0)$;
\item ${\Cal A}^\Sigma_\omega=\bigcup_{\sigma\in\Sigma}\, {\Cal K}^\sigma_\omega(0)$.
\endroster
\qed
\endproclaim

As in \cite{\rfa chepvish..}, in order to apply Theorem \rft ChepVish.., we need to check that
$P_\omega(t)(\sigma,u)$  is continuous, bounded, pointwise-dissipative and asymptotically
compact. Boundedness and pointwise-dissipativeness are a straightforward consequence of
Proposition \rft absorb1... For continuity and asymptotic compactness, we need some
preliminary lemmas.

\proclaim{Lemma \dft converg1..}
Let  $(\alpha_{ij}^1(\cdot))$ and $(\alpha_{ij}^2(\cdot))\in\Sigma_1$.
For $k=1,2$, let
$V^k_{\omega}(t,s)$ be the linear process in $L^2(\R^N)$ generated by the equation
$$
u_t=\sum_{i,j=1}^N \alpha_{ij}^k(\omega t)\partial_i\partial_j u.
$$
Then
\roster
\item there exists a continuous function $\rho\colon\R_+\to\R_+$, $\rho(q)\to 0$ as $q\to 0$,
such that, for any $\omega> 0$, for $u\in L^2(\R^N)$ and for $t>s$,
$$
\|V^1_{\omega}(t,s)u-V^2_\omega(t,s) u\|_{H^1}\leq \left(1+{{1}\over{(t-s)^{1/2}}}\right)
\rho(\|(a^1_{ij})-(a^2_{ij})\|_{\infty})\|u\|_{L^2};
$$
\item there exists a continuous function $\chi\colon\R_+\to\R_+$, $\chi(q)\to 0$ as $q\to 0$,
such that, for any $\omega> 0$, for $u\in H^1(\R^N)$ and for $t\geq s$,
$$
\|V^1_{\omega}(t,s)u-V^2_\omega(t,s) u\|_{H^1}\leq
\chi(\|(a^1_{ij})-(a^2_{ij})\|_{\infty})\|u\|_{H^1}.
$$
\endroster
\endproclaim
\demo{Proof}
Let $u\in L^2(\R^N)$. By \rff norm.. and \rff rappr.., for any $R>0$ we have
$$
\alignedat1
&\|V^1_\omega(t,s)u-V^2_\omega(t,s)u\|_{H^1}^2\\
&=\int_{\R^N}(1+|\xi|^2)
[\,\exp(-\int_s^t\sum_{i,j} \alpha_{ij}^1(\omega p)\xi_i\xi_j\,dp)\\
&\hbox to3cm{\null}- \exp(-\int_s^t\sum_{i,j}
\alpha_{ij}^2(\omega p)\xi_i\xi_j\,dp)\, ]^2 ({\Cal F}u)(\xi)^2\, d\xi\\
&\leq 2\int_{\{|\xi|\geq R\}}(1+|\xi|^2)\exp(-2\nu_0|\xi|^2(t-s))({\Cal F}u)(\xi)^2\,d\xi\\
&+\int_{\{|\xi|\leq R\}}(1+|\xi|^2)\exp(-2\nu_0|\xi|^2(t-s))\\
&\hbox to3cm{\null}[\,\exp(-\sum_{i,j}\int_s^t( \alpha_{ij}^1(\omega p)-\alpha_{ij}^2(\omega
p))\, dp\,\xi_i\xi_j)- 1\, ]^2  ({\Cal F}u)(\xi)^2\, d\xi\\
&=: S_1+S_2.
\endalignedat
$$
Choose $R:=k^{1/2}(t-s)^{-1/2}$, $k$ to be determined. If $k\geq1/(2\nu_0)$, we have
$$
S_1\leq 2\sup_{z\geq 2\nu_0
k}\left(1+{{z}\over{2\nu_0(t-s)}}\right)e^{-z}\|u\|_{L^2}^2
\leq 2\left(1+{1\over{(t-s)}}\right)ke^{-2\nu_0 k}\|u\|_{L^2}^2.
$$
On the other hand,
$$
S_2\leq \left(1+{{1}\over{(t-s)}}\right)k
[\,e^{\|(\alpha_{ij}^1)-(\alpha_{ij}^2)\|_{\infty}k} -1\,]^2\|u\|_{L^2}^2.
$$
Choosing
$$
k=\|(\alpha^1_{ij})-(\alpha^2_{ij})\|_\infty^{-1/3}
$$
we obtain the desired result with
$$
\rho(q)=\sqrt 2 q^{-1/6}e^{-\nu_0 q^{-1/3}}+ q^{-1/6}[e^{q^{2/3}}-1], \quad q\leq 8\nu_0^3.
$$
If $u\in H^1(\R^N)$ we argue in the same way: we get
$$
\|V^1_\omega(t,s)u-V^2_\omega(t,s)u\|_{H^1}^2\leq S_1+S_2,
$$
where
$$
S_1=2e^{-2\nu_0 k}\|u\|_{H^1}^2
$$
and
$$
S_2=[\,e^{\|(\alpha_{ij}^1)-(\alpha_{ij}^2)\|_{\infty}k} -1\,]^2\|u\|_{H^1}^2.
$$
Again, choosing $k$ as above,
we obtain the desired result.
\qed\enddemo

\proclaim{Lemma \dft contin1..}
Let $\sigma\in\Sigma$ and let $(\sigma_n)_{n\in\N}$ be a sequence in
$\Sigma$, such that $\sigma_n\to\sigma$ as $n\to\infty$. Let $(t_n)_{n\in\N}$ and
$(s_n)_{n\in\N}$ be two sequences of real numbers, with $t_n\geq s_n$ for all $n$ and
assume that $t_n\to t$ and $s_n\to s$ as $n\to\infty$. Finally, let $u\in H^1(\R^N)$ and let
$(u_n)_{n\in\N}$ be a  bounded sequence in $H^1(\R^N)$. Then, for any $\omega>0$,
\roster
\item if $u_n\to u$ in $L^2(\R^N)$ and $t>s$,
$$
\|\Pi^{\sigma_n}_\omega(t_n,s_n)u_n-\Pi^\sigma_\omega(t,s)u\|_{H^1}\to 0\quad\text{as
$n\to\infty$;}
$$
\item if $u_n\to u$ in $H^1(\R^N)$,
$$
\|\Pi^{\sigma_n}_\omega(t_n,s_n)u_n-\Pi^\sigma_\omega(t,s)u\|_{H^1}\to 0\quad\text{as
$n\to\infty$.}
$$
\endroster
\endproclaim
\demo{Proof}
First, let us notice that
$$
\Pi^{\sigma_n}_\omega(t_n,s_n)u_n=\Pi^{T_\omega(s_n-s)\sigma_n}_\omega(t_n-(s_n-s),s)u_n.
$$
Since $T_\omega(s_n-s)\sigma_n\to\sigma$ in $\Sigma$ and $t_n-(s_n-s)\to t$ as $n\to\infty$,
we can assume without loss of generality that $s_n=s$ for all $n$.
Let's write
$$
v_n(t):= \Pi_{\omega}^{\sigma_n}(t,s)u_n,
$$
$$
v(t):= \Pi_{\omega}^{\sigma}(t,s)u.
$$
We introduce the following notations:
$$
\alignedat1
\sigma_n(\tau)=:&((\alpha_{ij}^n (\tau)), \alpha_0^n(\tau), \varphi_n(\tau), \gamma_n(\tau)),\\
\sigma(\tau)=:&((\alpha_{ij}(\tau)), \alpha_0(\tau), \varphi(\tau), \gamma(\tau)),
\endalignedat
$$
and
$$
\alignedat1
A_n:=& \sup_{\tau \in \R} |(\alpha_{ij}^n(\tau))-(\alpha_{ij}(\tau))|,\\
B_n:=& \sup_{\tau \in \R} |\alpha_0^n(\tau) - \alpha_0(\tau)|,\\
D_n:=& \sup_{\tau \in \R} \sup_{u \in \R}
{{|(\varphi_n)_u(\tau,u)-\varphi_u(\tau,u)|}\over{1+|u|^{\beta}}},\\
E_n:=& \sup_{\tau \in \R} \|\gamma_n(\tau)-\gamma(\tau)\|_{L^2}.
\endalignedat
$$
Notice that $A_n$, $B_n$, $D_n$ and $E_n$ tend to zero as $n \rightarrow +\infty$.\par
\noindent Moreover, let's observe that for every $\tau \in \R$ we have
$$
\|\hat{\varphi_n}(\tau, u)-\hat{\varphi}(\tau,u)\|_{L^2}\leq D_n (\|u\|_{L^2}+
\|u\|_{H^1}^{\beta +1}).
$$
Finally, in view of Proposition \rft absorb1.., there exists $\tilde K>0$ such that, for every
$t\geq s$,
$$
\|v_n(t)\|_{H^1} \leq \tilde K \quad\text{for all
$n\in\N$ },
$$
$$
\|v(t)\|_{H^1} \leq \tilde K.
$$
Let $T$ be a positive number; for $0<t-s<T$ we have
$$ \multline
v_n(t)-v(t)= V^n_{\omega}(t,s)u_n - V_{\omega}(t,s)u  \\ +\int_s^t V^n_{\omega}(t,p) [ -
\alpha_0^n(\omega p) v_n(p) + \hat{\varphi}_n(\omega p, v_n(p))+ \gamma_n(\omega p)] \, dp  \\ -
\int_s^t V_{\omega} (t,p)[-\alpha(\omega p)v(p)+ \hat{\varphi}(\omega p,v(p))+ \gamma (\omega
p)]\,dp. \endmultline
$$
Hence
$$ \multline
\|v_n(t)-v(t)\|_{H^1} \leq \|V_\omega^n(t,s)[u_n-u]\|_{H^1}  \\
+\|[V_{\omega}^n(t,s)-V_n(t,s)]u\|_{H^1}+ I_1 +I_2 +I_3,
\endmultline
$$
where
$$
I_1:= \int_s^t \|(V_{\omega}^n(t,p)-V_{\omega}(t,p))[\alpha_0^n(\omega
p)v_n(p)+\hat{\varphi}_n(\omega p,v_n(p)) +\gamma_n(\omega p)]\|_{H^1}\,dp,
$$
$$ \multline
I_2:=\int_s^t \|V_{\omega}(t,p)[(-\alpha_0^n(\omega p)+\alpha_0(\omega p))v_n(p) \\+
\hat{\varphi}_n(\omega p,v_n(p))- \hat{\varphi}(\omega p,v_n(p))+ \gamma_n(\omega p)-\gamma
(\omega p)]\|_{H^1}\, dp,
\endmultline
$$
$$
I_3:= \int_s^t \|V_{\omega}(t,p)[-\alpha_0(\omega p)(v_n(p)-v(p))+
\hat{\varphi}(\omega p,v_n(p))- \hat{\varphi}(\omega p,v(p))]\|_{H^1}\,dp.
$$
First of all, let's observe that, thanks to Lemma \rft converg1..,
$$
\|[V_{\omega}^n(p)- V_{\omega}(p)]u \|_{H^1} \leq \chi(A_n) \|u\|_{H^1}\leq \chi(A_n)\tilde K.
$$
As for $I_1$, Lemma \rft converg1.. implies that
$$
\multline
I_1 \leq \int_s^t (1+(t-p)^{-1/2}) \rho (A_n) \|\alpha_0^n(\omega
p)v_n(p)+\hat{\varphi_n}(\omega p,v_n(p)) +\gamma_n(\omega p)\|_{L^2} \, dp \\
\leq \int_s^t (1+(t-p)^{-1/2}) \rho (A_n) (C \tilde K+ \tilde C (\tilde K
+\tilde{K}^{\beta +1})+C)\, dp \leq Q_1 \rho (A_n),
\endmultline
$$
where $Q_1$ is a positive constant depending on $T$. Analogously, \rff expest3.. implies that
$$
\multline
I_2 \leq \int_s^t M (1+(t-p)^{-1/2}) \|(-\alpha_0^n(\omega p)+\alpha_0(\omega p))v_n(p) \\+
\hat{\varphi}_n(\omega p,v_n(p))- \hat{\varphi}(\omega p,v_n(p))+  \gamma_n(\omega p)-\gamma
(\omega p) \|_{L^2} \,dp \\
\leq \int_s^t M (1+(t-p)^{-1/2}) (B_n \tilde K +D_n (\tilde K +\tilde{K}^{\beta
+1})+E_n) \, dp \leq Q_2(B_n+D_n+E_n),
\endmultline
$$
where $Q_2$ is a positive constant depending on $T$.
Finally, \rff expest3.. implies
$$\multline
I_3 \leq \int_s^t M (1+(t-p)^{-1/2}) \|-\alpha_0(\omega p)(v_n(p)-v(p))\\
+\hat{\varphi}(\omega p,v_n(p))- \hat{\varphi}(\omega p,v(p)) \|_{L^2} \,dp \\
\leq \int_s^t M (1+(t-p)^{-1/2}) [ C\|v_n(p)-v(p)\|_{H^1}+ \tilde C
(1+2\tilde{K}^{\beta})\|v_n(p)-v(p)\|_{H^1}] \, dp \\
\leq Q_3 \int_s^t (t-p)^{-1/2} \|v_n(p)-v(p)\|_{H^1}\,dp,
\endmultline
$$
where $Q_3$ is a positive constant depending on $T$.\par
As a consequence,
$$
\|v_n(t)-v(t)\|_{H^1}\leq \|V^n_{\omega}(t,s)(u_n-u)\|_{H^1}+ Z_n +Q_3 \int_s^t
(t-p)^{-1/2} \|v_n(p)-v(p)\|_{H^1} \, dp,
$$
where $Z_n \rightarrow 0$ as $n\rightarrow +\infty$.\par
In case (1) we have, due to \rff expest3..,
$$
\| V^n_{\omega}(t,s)(u_n-u)\|_{H^1}\leq M(1+(t-s)^{-1/2}) \|u_n-u\|_{L^2}\leq \tilde M
(t-s)^{-1/2}\|u_n-u\|_{L^2},
$$
hence
$$
\|v_n(t)-v(t)\|_{H^1} \leq (t-s)^{-1/2} F_n +Q_3 \int_s^t (t-p)^{-1/2}
\|v_n(p)-v(p)\|_{H^1}\,dp,
$$
where $F_n \rightarrow 0$ as $n\rightarrow +\infty$, and by the singular version of Gronwall's
inequality (see \cite{\rfa He.., Th. 7.1.1})
$$
\|v_n(t)-v(t)\|_{H^1}\leq  Q F_n (t-s)^{-1/2},
$$
where $Q$ is a positive constant. This implies $v_n(t) \rightarrow v(t)$ in $H^1$
uniformly on $[s+\delta,s+T]$ for every $T>\delta >0$, and proves (1).\par
In case (2), \rff expest2.. implies that
$$
\|V_{\omega}^n(t,s)(u_n-u)\|_{H^1}\leq M\|u_n-u\|_{H^1},
$$
hence
$$
\|v_n(t)- v(t)\|_{H^1}\leq \tilde{F_n} + \tilde{Q}_3 \int_s^t (t-p)^{-1/2}
\|v_n(p)-v(p)\|_{H^1} \,dp,$$
where $\tilde{Q}_3$ is a positive constant and $\tilde{F_n} \rightarrow 0$ as $n\rightarrow
+\infty$. Again by the singular version of Gronwall's inequality
$$
\|v_n(t)-v(t)\|_{H^1} \leq \tilde Q \tilde{F_n},
$$
where $\tilde Q$ is a positive constant. This implies $v_n(t) \rightarrow v(t)$ in $H^1$
uniformly on $[s,s+T]$, and  proves (2).
\qed\enddemo

We recall the following
\proclaim{Definition \dft defascomp..}
A bounded semigroup $P(t)$ acting on a complete metric space $X$ is said to be asymptotically
compact iff for every bounded sequence $(u_n)_{n\in\N}$ and for every sequence
$(t_n)_{n\in\N}$,
$t_n\to+\infty$ as $n\to\infty$, there exists $u_\infty\in X$ such that, up to a subsequence,
$P(t_n)u_n\to u_\infty$ as $n\to\infty$.
\endproclaim

Now we can prove
\proclaim{Proposition \dft ascomp1..}
The semigroup $P_\omega(t)$ is continuous and asymptotically compact on $\Sigma\times
H^1(\R^N)$.
\endproclaim
\demo{Proof}
The continuity of $P_\omega(t)$ is a straightforward consequence of Lemma \rft contin1.. and we
omit the easy proof.

In order to prove the asymptotic compactness of $P_\omega(t)$, we take a bounded
sequence $((\sigma_n,u_n))_{n\in\N}$ in $\Sigma\times H^1(\R^N)$. Let $R>0$ be such that
$\|u_n\|_{H^1}\leq R$ for all $n\in\N$. We seek for $(\sigma_\infty,u_\infty)\in\Sigma\times
H^1(\R^N)$ such that, up to a subsequence,
$$
P_\omega(t_n)(\sigma_n,u_n)=(T_\omega(t_n)\sigma_n,\Pi^{\sigma_n}_\omega(t_n,0)u_n)\to
(\sigma_\infty,u_\infty)\quad\text{in $\Sigma\times H^1(\R^N)$ }
$$
as $n\to\infty$. First of all, since $\Sigma$ is compact, we can assume, without loss of
generality, that there exists $\cl\sigma_\infty\in\Sigma$ such that
$T_\omega(t_n-1)\sigma_n\to\cl\sigma_\infty$ and
$T_\omega(t_n)\sigma_n\to T_\omega(1)\cl\sigma_\infty=:\sigma_\infty$ as $n\to\infty$.
Moreover, since the sequence $(u_n)_{n\in\N}$ is bounded in $H^1(\R^N)$, by Proposition
\rft absorb1.. the set
$$
\{\,\Pi^{\sigma_n}_\omega(t_n,0)u_n\mid n\in\N\,\}
$$
is bounded, and
hence weakly compact in $H^1(\R^N)$. So, passing to a subsequence if necessary, we can assume
that there exists
$u_\infty\in H^1(\R^N)$ such that
$$
\Pi^{\sigma_n}_\omega(t_n,0)u_n\weakto u_\infty \quad\text{in $H^1(\R^N)$ as $n\to\infty$.}
$$
We must show that the convergence is actually strong in $H^1(\R^N)$.\par
We claim first that $\Pi^{\sigma_n}_\omega(t_n,0)u_n\to u_\infty $ in the strong $L^2$-topology.
To this end, it is enough to show that the set
$$
\{\,\Pi^{\sigma_n}_\omega(t_n,0)u_n\mid n\in\N\,\}
$$
is relatively compact in the strong $L^2$ topology, or equivalently that it is
totally bounded. This is a consequence of Lemma \rft ball1.. and of Rellich theorem.
Let $\eta>0$; by Lemma \rft ball1.., there exists $k>0$ and $\cl n\in\N$, depending on
$R$ and $\eta$, such
that
$$
\int_{\{|x|>k\}}|\Pi^{\sigma_n}_\omega(t_n,0)u_n(x)|^2 \,dx\leq\eta\quad\text{for all $n\geq
\cl n$.}
$$
 We introduce the operator $\Xi\co L^2(\R^N)\to L^2(\R^N)$,
$$
(\Xi u)(x):=\cases u(x)&\text{if $|x|\leq
k$}\\0&\text{if $|x|>k$}\endcases
$$
Then we have
$$
\multline
\{\,\Pi^{\sigma_n}_\omega(t_n,0)u_n\mid n\in\N\,\}
=\{\,\Xi\,\Pi^{\sigma_n}_\omega(t_n,0)u_n+(I-\Xi)\Pi^{\sigma_n}_\omega(t_n,0)u_n\mid n\in\N\,\}\\
\subset
\{\,\Xi\,\Pi^{\sigma_n}_\omega(t_n,0)u_n\mid
n\in\N\,\}+\{\,(I-\Xi)\Pi^{\sigma_n}_\omega(t_n,0)u_n\mid n\in\N\,\}\\
\subset B_\eta(0)+\{\,\Xi\,\Pi^{\sigma_n}_\omega(t_n,0)u_n\mid
n\in\N\,\}\endmultline
$$
where $B_\eta(0)$ is the ball of radius $\eta$ centered at $0$ in $L^2(\R^N)$.
The set
$$
\{\,\Xi\,\Pi^{\sigma_n}_\omega(t_n,0)u_n\mid n\in\N\,\}
$$
consists of functions of $L^2(\R^N)$ which are equal to zero outside the ball of radius $k$ in
$\R^{N}$ and whose restriction to the same ball is in $H^1$. On the other hand, the $H^1$-norm of
these functions is uniformly bounded. Then, by Rellich Theorem, we deduce that the set
$\{\,\Xi\,\Pi^{\sigma_n}_\omega(t_n,0)u_n\mid n\in\N\,\}$ is precompact in $L^2(\R^N)$. Hence we
can cover it by a finite number of balls of radius $\eta$ in $L^2(\R^N)$. This implies that the
set $\{\,\Pi^{\sigma_n}_\omega(t_n,0)u_n\mid n\in\N\,\}$ is totally bounded and hence
precompact in $L^2(\R^N)$. The claim is proved.\par
The same conclusions obviously hold also for the set
$$
\{\,\Pi^{\sigma_n}_\omega(t_n-1,0)u_n\mid n\in\N\,\};
$$
so there exists $\cl u_\infty\in H^1(\R^N)$ such that, up to a subsequence,
$$
\Pi^{\sigma_n}_\omega(t_n-1,0)u_n\to \cl u_\infty \quad\text{in $L^2(\R^N)$ as $n\to\infty$.}
$$
Finally, by Lemma \rft contin1.., we have
$$
\multline
\Pi^{\sigma_n}_\omega(t_n,0)u_n
=\Pi^{\sigma_n}_\omega(t_n,t_n-1)\Pi^{\sigma_n}_\omega(t_n-1,0)u_n\\
=\Pi^{T_\omega(t_n-1)\sigma_n}_\omega(1,0)\Pi^{\sigma_n}_\omega(t_n-1,0)u_n
\to \Pi^{\cl\sigma_\infty}_\omega(1,0)\cl u_\infty\quad\text{in $H^1(\R^N)$ as $n\to\infty$.}
\endmultline$$
It follows that $u_\infty=\Pi^{\cl\sigma_\infty}_\omega(1,0)\cl u_\infty$ and
$\Pi^{\sigma_n}_\omega(t_n,0)u_n\to u_\infty$ in $H^1(\R^N)$ as $n\to\infty$. The proof is
complete.
\qed\enddemo

Finally, combining Theorem \rft ChepVish.. and Proposition \rft ascomp1..,  we have:
\proclaim{Theorem \dft exattr..}
The family of processes $\{\,\Pi^\sigma_\omega\mid\sigma\in\Sigma\,\}$ in $H^1(\R^N)$
possesses a compact $\Sigma$-uniform attractor ${\Cal A}^\Sigma_\omega$. As a point set,
$$
{\Cal A}^\Sigma_\omega=\bigcup_{\sigma\in\Sigma}\, {\Cal K}^\sigma_\omega(0),
$$
where ${\Cal K}^\sigma_\omega(0)$ is the kernel section introduced in Definition \rft kernel...
In other words, ${\Cal A}^\Sigma_\omega$ is the union of all the full bounded trajectories of
$\Pi^\sigma_\omega$, $\sigma\in\Sigma$.
\qed
\endproclaim

\newsection\head \the\secnumber. Behaviour as $\omega\to+\infty$\endhead
In this section we shall investigate the behaviour of the solutions of \rff equazione.. as
$\omega\to+\infty$. As we explained in the Introduction, we expect that the averaged equation
$$
u_t=\sum_{i,j=1}^N \cl a_{ij}\partial_i\partial_j u-\cl a_0 u+\cl f(u)+\cl g(x)
\tag\dff averaged..
$$
behaves like a `limit' equation of \rff equazione... Roughly speaking, this means that the
solutions of \rff equazione.. with initial datum $u_0\in H^1(\R^N)$, as $\omega\to+\infty$,
converge in some sense to the solution of
\rff averaged.. with the same initial datum. Moreover, we claim that the attractor of \rff
equazione.. is $H^1$-close to that of \rff averaged.. for sufficiently large $\omega$.

Let us denote by $\cl A\colon H^2(\R^N)\to L^2(\R^N)$ the self-adjoint positive operator defined
by
$$
\cl A u:=-\sum_{i,j=1}^N\cl a_{ij}\partial_i\partial_j u,\quad u\in H^2(\R^N);
$$
we denote by $e^{-\cl At}$ the analytic semigroup generated by $\cl A$. Then equation
\rff averaged.. can be written as an abstract parabolic equation in $L^2(\R^N)$, namely
$$
\dot u=-\cl Au-\cl a_0u+\Hat{\Bar f}(u)+\cl g.\tag\dff limeq..
$$
Equation \rff limeq.. defines a global semiflow $\pi$ in $H^1(\R^N)$: in fact, as we already
observed in Section 3, all a-priori estimates of Section 2 are independent of $\omega$ and
$\sigma\in\Sigma$, and are valid also for the averaged equation \rff averaged...
So the semiflow $\pi$ possesses a compact global attractor $\Cal A$.

We begin with a convergence result for the linear problems associated to \rff equazione.. and
\rff averaged..:

\proclaim{Proposition \dft linaver..}
For $(\alpha_{ij}(\cdot))\in \Sigma_1$ and $\omega>0$, let $V^\alpha_\omega(t,s)$ be the linear
process  in $L^2(\R^N)$ generated by the equation
$$
u_t=\sum_{i,j=1}^N \alpha_{ij}(\omega t)\partial_i\partial_j u.
$$
Moreover, let $e^{-\bar A t}$ be the linear semigroup in $L^2(\R^N)$ generated by the equation
$$
u_t=\sum_{i,j=1}^N \bar a_{ij}\partial_i\partial_j u.
$$

There exists a bounded, continuous and decreasing function $\theta\colon \R_+\to\R_+$,
$\theta(q)\to 0$ as $q\to\infty$, such that, for $u\in L^2(\R^N)$ and for $t>s$,
$$
\|V^\alpha_\omega(t,s)u-e^{-\bar A(t-s)}u\|_{H_1}\leq \left(1+{1\over{(t-s)^{1/2}}}\right)
\theta(\omega(t-s))\|u\|_{L^2}
$$
for any $(\alpha_{ij}(\cdot))\in \Sigma_1$ and $\omega>0$.
\endproclaim
\demo{Proof}
Let $u\in L^2(\R^N)$ and $t>s$. By \rff norm.. and \rff rappr.., for  any $R>0$ we have
$$
\alignedat1
&\|V^\alpha_{\omega}(t,s)u-e^{-\bar A(t-s)}u\|_{H^1}^2\\
&=\int_{\R^N}(1+|\xi|^2)
[\,\exp(-\int_s^t\sum_{i,j} \alpha_{ij}(\omega p)\xi_i\xi_j\,dp)\\
&\hbox to3cm{\null}- \exp(-\sum_{i,j}
\bar a_{ij}\xi_i\xi_j(t-s))\, ]^2 ({\Cal F}u)(\xi)^2\, d\xi\\
&\leq 2\int_{\{|\xi|\geq R\}}(1+|\xi|^2)\exp(-2\nu_0|\xi|^2(t-s))({\Cal F}u)(\xi)^2\,d\xi\\
&+\int_{\{|\xi|\leq R\}}(1+|\xi|^2)\exp(-2\nu_0|\xi|^2(t-s))\\
&\hbox to3cm{\null}[\,\exp(-\sum_{i,j}\int_s^t( \alpha_{ij}(\omega p)-\bar a_{ij})\,
dp\,\xi_i\xi_j)- 1\, ]^2  ({\Cal F}u)(\xi)^2\, d\xi\\
&=: S_1+S_2.
\endalignedat
$$
Choose $R:=k^{1/2}(t-s)^{-1/2}$, $k$ to be determined. If $k\geq 1/2\nu_0$,
$$
S_1\leq 2\sup_{z\geq k}\left(1+{{z}\over{(t-s)}}\right)e^{-2\nu_0z}\|u\|^2_{L^2}
\leq 2\left(1+{{1}\over{(t-s)}}\right)ke^{-2\nu_0k}\|u\|^2_{L^2}.
$$
In order to estimate $S_2$, we observe that
$$
\multline
|\int_s^t( \alpha_{ij}(\omega p)-\bar a_{ij})\, dp|\\
=(t-s)\left|{1\over{\omega(t-s)}}\int_{\omega s}^{\omega s+\omega(t-s)}
(\alpha_{ij}(p)-\bar a_{ij})\, dp\right|
\leq (t-s)\mu(\omega(t-s)).
\endmultline
$$
Then
$$
\multline
S_2\leq \left(1+{{k}\over{(t-s)}}\right)
\int_{\{|\xi|^2\leq k(t-s)^{-1}\}}
[\,e^{Nk\mu(\omega(t-s))}- 1\, ]^2  ({\Cal F}u)(\xi)^2\, d\xi\\
\leq \left(1+{{1}\over{(t-s)}}\right)k
[\,e^{Nk\mu(\omega(t-s))}- 1\, ]^2  \|u\|_{L^2}^2
\endmultline
$$
By the mean value theorem, we get
$$
\multline
S_2\leq \left(1+{{1}\over{(t-s)}}\right)k
[\,e^{Nk\mu(\omega(t-s))}Nk\mu(\omega(t-s))\, ]^2  \|u\|_{L^2}^2\\
=\left(1+{{1}\over{(t-s)}}\right)N^2k^3
e^{2Nk\mu(\omega(t-s))}\mu(\omega(t-s))^2  \|u\|_{L^2}^2
\endmultline
$$
Now, set $\mu_\infty:=\sup_{q\geq 0}\mu(q)$ and take
$$
k:={1\over{2\nu_0}}\left({{\mu_\infty}\over{\mu(\omega(t-s))}}\right)^{1/2}
$$
With this choice of $k$, we obtain that there exist positive constants $\cl\kappa$, $\mu_1$ and
$\mu_2$, depending only on $N$, $\mu_\infty$ and $\nu_0$, such that
$$
\multline
S_1+S_2\leq \cl\kappa\left(1+{1\over{(t-s)}}\right)\\
\left(\mu(\omega(t-s))^{-1/2}e^{-\mu_1\mu(\omega(t-s))^{-1/2}}+
\mu(\omega(t-s))^{1/2}e^{\mu_2\mu(\omega(t-s))^{1/2}}\right)
\|u\|^2_{L^2}
\endmultline
$$
The conclusion follows if we define
$$
\theta(q):=\cl\kappa^{1/2}\left(\mu(q)^{-1/2}e^{-\mu_1\mu(q)^{-1/2}}+
\mu(q)^{1/2}e^{\mu_2\mu(q)^{1/2}}\right)^{1/2}
$$
and, if necessary, modify it on some bounded interval in order to make it decreasing on $\R_+$.
The proof is complete.
\qed\enddemo

\proclaim{Corollary \dft converg2..}
Let $(\omega_n)_{n\in\N}$ be a sequence of positive numbers, $\omega_n\to+\infty$ as
$n\to\infty$. Let $(\alpha_{ij}^n(\cdot))_{n\in\N}$ be a sequence in $\Sigma_1$ and let
$V^n_{\omega_n}(t,s)$ be the linear process in $L^2(\R^N)$ generated by the equation
$$
u_t=\sum_{i,j=1}^N \alpha_{ij}^n(\omega_n t)\partial_i\partial_j u.
$$
Fix $0<\delta<T$ and
take a sequence $(u_n)_{n\in\N}$ in $L^2(\R^N)$, $u_n\to u$ in $L^2(\R^N)$. Then
$$
\sup_{s\in\R}\,\sup_{t\in[s+\delta,s+T]}\|V^n_{\omega_n}(t,s)u_n-e^{-\bar A(t-s)}u\|_{H^1}\to 0
\quad\text{as $n\to\infty$.}
$$
\endproclaim
\demo{Proof}
We have
$$
\multline
\|V^n_{\omega_n}(t,s)u_n-e^{-\bar A(t-s)}u\|_{H^1}\\
\leq \|V^n_{\omega_n}(t,s)u_n-V^n_{\omega_n}(t,s)u\|_{H^1}+
\|V^n_{\omega_n}(t,s)u-e^{-\bar A(t-s)}u\|_{H^1}\\
\leq M(1+(t-s)^{-1/2})\left(\|u_n-u\|_{L^2}+\theta(\mu(\omega_n(t-s)))\|u\|_{L^2}\right)\\
\leq M(1+\delta^{-1/2})\left(\|u_n-u\|_{L^2}+\theta(\mu(\omega_n\delta))\|u\|_{L^2}\right)\to 0
\quad\text{as $n\to\infty$,}
\endmultline
$$
and the corollary is proved.\qed\enddemo

If we deal with a fixed $u\in H^1(\R^N)$, we obtain uniform convergence on the whole interval
$[s,s+T]$. Indeed, we have the following

\proclaim{Proposition \dft linaver2..}
Let $(\omega_n)_{n\in\N}$ be a sequence of positive numbers, $\omega_n\to\infty$ as $n\to\infty$.
Let $(\alpha_{ij}^n(\cdot))_{n\in\N}$ be a sequence in $\Sigma_1$ and let
$V^n_{\omega_n}(t,s)$ be the linear process in $L^2(\R^N)$ generated by the equation
$$
u_t=\sum_{i,j=1}^N \alpha_{ij}^n(\omega_n t)\partial_i\partial_j u.
$$
Finally, let $u\in H^1(\R^N)$. Then, for any $T>0$,
$$
\sup_{s\in\R}\,\sup_{t\in[s,s+T]}\|V^n_{\omega_n}(t,s)u-e^{-\bar A(t-s)}u\|_{H^1}\to 0
\quad\text{as $n\to\infty$.}
$$
\endproclaim
\demo{Proof}
Let $u\in H^1(\R^N)$ and $t>s$. Arguing like in the proof of Proposition \rft linaver.., for
any $R>0$ we have
$$
\alignedat1
&\|V^n_{\omega_n}(t,s)u-e^{-\bar A(t-s)}u\|_{H^1}^2
\leq 2\int_{\{|\xi|\geq R\}}(1+|\xi|^2)({\Cal F}u)(\xi)^2\,d\xi\\
&+\int_{\{|\xi|\leq R\}}(1+|\xi|^2)[\,\exp(-\sum_{i,j}\int_s^t
(\alpha_{ij}^n(\omega_n p)-\bar a_{ij})\, dp\,\xi_i\xi_j)- 1\, ]^2  ({\Cal F}u)(\xi)^2\, d\xi.\\
\endalignedat
$$
Since
$$
|\int_s^t( \alpha_{ij}(\omega p)-\bar a_{ij})\, dp|
\leq (t-s)\mu(\omega(t-s)),
$$
we obtain
$$
\multline
\|V^n_{\omega_n}(t,s)u-e^{-\bar A(t-s)}u\|_{H^1}^2\\
\leq 2\int_{\{|\xi|\geq R\}}(1+|\xi|^2)({\Cal F}u)(\xi)^2\,d\xi+
\left(e^{N(t-s)\mu(\omega_n(t-s))R^2}-1\right)^2\|u\|^2_{H^1}
\endmultline
$$
Now, given $\epsilon>0$, we choose $R>0$ (depending on $u$ and $\epsilon$) such that
$$
\int_{\{|\xi|\geq R\}}(1+|\xi|^2)({\Cal F}u)(\xi)^2\,d\xi\leq \epsilon.
$$
Let $\delta$ be a positive number, depending on $R$, $\|u\|_{H^1}$ and $\epsilon$,
such that $\delta<T$ and
$$
\left(e^{N\delta\mu_\infty R^2}-1\right)^2\|u\|_{H^1}^2\leq\epsilon.
$$
Then, for $t-s<\delta$,
$$
\left(e^{N(t-s)\mu(\omega_n(t-s))R^2}-1\right)^2\|u\|_{H^1}^2
\leq \left(e^{N\delta\mu_\infty R^2}-1\right)^2\|u\|_{H^1}^2\leq\epsilon.
$$
On the other hand, if $\delta\leq (t-s)\leq T$, we have
$$
\left(e^{N(t-s)\mu(\omega_n(t-s))R^2}-1\right)^2
\leq \left(e^{NT\mu(\omega_n\delta)R^2}-1\right)^2.
$$
As a consequence, given $\epsilon>0$, we can find $R$ and $\delta$ (depending
on $\epsilon$) such that, for all $n\in\N$,
$$
\sup_{s\in\R}\sup_{t\in[s,s+T]}\|V^n_{\omega_n}(t,s)u-e^{-\bar A(t-s)}u\|_{H^1}^2\leq
3\eps+\left(e^{NT\mu(\omega_n\delta)R^2}-1\right)^2\|u\|^2_{H^1}.
$$
The conclusion follows by letting $n\to\infty$.
\qed\enddemo

\remark{Remark}
The convergence in Proposition \rft linaver2.. is not uniform with respect to $u$ in a bounded
subset of
$H^1(\R^N)$. As a matter of fact, if we try to repeat the arguments of Proposition \rft
linaver.., we see that there exists a bounded, continuous and decreasing function
$\theta\colon
\R_+\to\R_+$,
$\theta(q)\to 0$ as $q\to\infty$, such that, for $u\in H^1(\R^N)$ and for $t>s$,
$$
\|V^\alpha_\omega(t,s)u-e^{-\bar A(t-s)}u\|_{H_1}\leq
\theta(\omega(t-s))\|u\|_{H^1}
$$
for any $(\alpha_{ij}(\cdot))\in \Sigma_1$ and $\omega>0$. It is clear that this is not enough to
detect uniform convergence up to $t=s$, since there is still an initial layer one cannot get rid
of. This is due to the microlocal effect of the rapid oscillations of the coefficients
$a_{ij}(\omega p)$.
\endremark

\medskip

Now we can state our first `local' averaging result for the nonlinear equation \rff equazione..:

\proclaim{Theorem \dft contin2..}
Let $(\sigma_n)_{n\in\N}$ be a sequence in
$\Sigma$. Let $(t_n)_{n\in\N}$ and
$(s_n)_{n\in\N}$ be two sequences of real numbers, with $t_n> s_n$ for all $n$ and
assume that $t_n\to t$ and $s_n\to s$ as $n\to\infty$, with $t>s$. Let $u\in H^1(\R^N)$
and let
$(u_n)_{n\in\N}$ be a  bounded sequence in $H^1(\R^N)$ and assume that $u_n\to u$ in $L^2(\R^N)$.
Finally let $(\omega_n)_{n\in\N}$ be a sequence of positive numbers, $\omega_n\to+\infty$ as
$n\to\infty$. Then
$$
\|\Pi^{\sigma_n}_{\omega_n}(t_n,s_n)u_n-\pi(t-s)u\|_{H^1}\to 0\quad\text{as
$n\to\infty$.}
$$
\endproclaim

In order to prove Theorem \rft contin2.., we need the following
\proclaim{Lemma \dft tecnico1..}
Let $(\alpha_0^n,\phi_n,\gamma_n)_{n\in\N}$ be a sequence in
$\Sigma_2\times\Sigma_3\times\Sigma_4$. Let $(\omega_n)_{n\in\N}$ be a sequence of positive
numbers, $\omega_n\to+\infty$ as
$n\to\infty$. Finally, let $v\colon [s,s+T]\to H^1(\R^N)$ be a continuous function.
For $t\in[s,s+T]$ set
$$
G^1_n(t):=\int_s^t e^{-\bar A(t-p)}\left[\alpha_0^n(\omega_n p)-\cl a_0\right]v(p)\, dp;
$$
$$
G^2_n(t):=\int_s^t e^{-\bar A(t-p)}\left[\hat\phi_n(\omega_n p,v(p))-\Hat{\Bar f}(v(p))\right]\,
dp;
$$
$$
G^3_n(t):=\int_s^t e^{-\bar A(t-p)}\left[\gamma_n(\omega_n p)-\cl g\right]\, dp.
$$
Then $G^j_n(t)\to 0$ in $H^1(\R^N)$ uniformly on $[s,s+T]$ for $j=1,2,3$.
\endproclaim
\demo{Proof}
The proof of this lemma is essentially contained in \cite{\rfa il.., Th. 1.1} and, in a more
general setting, in \cite{\rfa He.., Th. 3.4.7}. We give
the details for sake of completeness. \par We start by considering $G^3_n$. First of all,
we observe that, for
$s<p<t$,
$$
\multline
{d\over{dp}}\left(e^{-\bar A(t-p)}\int_p^t\left[\gamma_n(\omega_n q)-\cl g\right]\,dq\right)\\
=\bar Ae^{-\bar A(t-p)}\int_p^t\left[\gamma_n(\omega_n q)-\cl g\right]\,dq+
e^{-\cl A(t-p)}\left[\gamma_n(\omega_n p)-\cl g\right].
\endmultline
$$
Since
$$
\multline
\|\bar Ae^{-\bar A(t-p)}\int_p^t\left[\gamma_n(\omega_n q)-\cl g\right]\,dq\|_{H^1}\\
\leq M(t-p)^{-3/2}(t-p)\bigl\|(\omega_n(t-p))^{-1}\int_{\omega_n
p}^{\omega_np+\omega_n(t-p)}\left[\gamma_n( q)-\cl g\right]\,dq\bigr\|_{L^2}\\
\leq M(t-p)^{-1/2}\mu(\omega_n(t-p))\in L^1(]\,p,t\,[)
\endmultline\tag\dff 27..
$$
and
$$
\|e^{-\cl A(t-p)}\left[\gamma_n(\omega_n p)-\cl g\right]\|_{H^1}\leq 2C(t-p)^{-1/2}\in
L^1(]\,p,t\,[),
$$
the `integration-by-part' formula
$$
\multline
\int_s^t e^{-\bar A(t-p)}\left[\gamma_n(\omega_n p)-\cl g\right]\, dp\\
=-e^{-\bar A(t-s)}\int_s^t\left[\gamma_n(\omega_n p)-\cl g\right]\,dp-
\int_s^t\bar Ae^{-\bar A(t-p)}\int_p^t\left[\gamma_n(\omega_n q)-\cl g\right]\,dq \, dp
\endmultline
$$
is valid. In view of \rff 27.., we get
$$
\multline
\|G^3_n(t)\|_{H^1}\leq
\left(1+{1\over(t-s)^{1/2}}\right)(t-s)
\mu(\omega_n(t-s))
+\int_s^t(t-p)^{-1/2}\mu(\omega_n(t-p))\,dp\\
\leq (t-s)^{1/2}\mu(\omega_n(t-s))+\int_s^t(t-p)^{-1/2}\mu(\omega_n(t-p))\,dp.
\endmultline
$$
Now let $\epsilon>0$. If $t-s\leq\epsilon$, a simple integration yields
$$
\|G^3_n(t)\|_{H^1}\leq 3\mu_\infty\epsilon^{1/2}.
$$
If $t-s\geq\epsilon$, we have
$$
\multline
\|G^3_n(t)\|_{H^1}\leq
T^{1/2}\mu(\omega_n\epsilon)+\mu(\omega_n\epsilon)\int_s^{t-\epsilon}(t-p)^{-1/2}\,dp
+\mu_\infty\int_{t-\epsilon}^t(t-p)^{-1/2}\,dp\\
\leq 3T^{1/2}\mu(\omega_n\epsilon)+2\mu_\infty\epsilon^{1/2}.
\endmultline
$$
It follows that, for all $n\in\N$,
$$
\sup_{t\in[s,s+T]}\|G^3_n(t)\|_{H^1}\leq 3\mu_\infty \eps^{1/2}+3T^{1/2}\mu(\omega_n\epsilon).
$$
The conclusion follows by letting $n\to\infty$.\par

Next we consider $G^1_n$. We assume first that $v(t)\equiv \cl v\in H^1(\R^N)$. Then, arguing
as above, we see that
$$
\multline
G^1_n(t)=-e^{-\bar A(t-s)}\int_s^t\left[\alpha^n_0(\omega_n p)-\cl a_0\right]\cl v\,dp\\
-\int_s^t\bar Ae^{-\bar A(t-p)}\int_p^t\left[\alpha^n_0(\omega_n q)-\cl a_0\right]\cl v\,dq 
\, dp,
\endmultline
$$
hence
$$
\|G^1_n(t)\|_{H^1}\leq
(t-s)\mu(\omega_n(t-s))\|\cl v\|_{H^1}
+\int_s^t\mu(\omega_n(t-p))\,dp\|\cl v\|_{H^1}.
$$
The same argument used for estimating $G^3_n$ shows that $G^1_n(t)\to 0$ in $H^1(\R^N)$
uniformly on $[s,s+T]$. One can easily see that the same is true if $v(t)$ is an
arbitrary bounded step function. The conclusion then follows by a density argument.\par

Finally, we consider $G^2_n$. Again we assume first that $v(t)\equiv\cl v\in H^1(\R^N)$. Then
we have
$$
\multline
G^2_n(t)=-e^{-\bar A(t-s)}\int_s^t\left[\hat\phi(\omega_n p,\cl v)-\Hat{\Bar f}(\cl
v)\right]\,dp\\ -\int_s^t\bar Ae^{-\bar A(t-p)}\int_p^t\left[\hat\phi(\omega_n q,\cl
v)-\Hat{\Bar f}(\cl v)\right]\,dq\, dp.
\endmultline
$$
By \rff 31.., we obtain
$$
\multline
\|G^2_n(t)\|_{H^1}\leq
K\left(1+{1\over(t-s)^{1/2}}\right)(t-s)
\mu(\omega_n(t-s))(\|\cl v\|_{L^2}+\|\cl v\|_{H^1}^{\beta+1})\\
+K\int_s^t(t-p)^{-1/2}\mu(\omega_n(t-p))(\|\cl v\|_{L^2}+\|\cl v\|_{H^1}^{\beta+1}) \,dp\\
\leq K\left( (t-s)^{1/2}\mu(\omega_n(t-s))+\int_s^t(t-p)^{-1/2}\mu(\omega_n(t-p))\,dp\right)
(\|\cl v\|_{L^2}+\|\cl v\|_{H^1}^{\beta+1}).
\endmultline
$$
Arguing as before, $G^2_n(t)\to 0$ in $H^1(\R^N)$ uniformly on $[s,s+T]$, and the same is
true if
$v(t)$ is an arbitrary bounded step function. The conclusion then follows again by a density
argument.
\qed\enddemo

\demo{Proof of Theorem \rft contin2..}
First, let us notice that
$$
\Pi^{\sigma_n}_{\omega_n}(t_n,s_n)u_n=\Pi^{T_{\omega_n}(s_n-s)\sigma_n}_{\omega_n}
(t_n-(s_n-s),s)u_n.
$$
Since $t_n-(s_n-s)\to t$ as $n\to\infty$,
we can assume without loss of generality that $s_n=s$ for all $n$.

Let's write
$$
v_n(t):=\Pi^{\sigma_n}_{\omega_n}(t,s)u_n,
$$
$$
v(t):=\pi(t-s)u.
$$
We recall that, in view of Proposition \rft absorb1.., there exists $\tilde K>0$ such that, for every
$t\geq s$,
$$
\|v_n(t)\|_{H^1} \leq \tilde K \quad\text{for all
$n\in\N$ },
$$
$$
\|v(t)\|_{H^1} \leq \tilde K.
$$
By the variation of constant formula we get
$$
\|v_n(t)-v(t)\|_{H^1} \leq \| V^n_{\omega_n}(t,s)u_n -e^{- \bar
A(t-s)}u\|_{H^1}+I_1(t)+I_2(t)+I_3(t),
$$
where
$$
\multline
I_1(t):= \int_s^t \|(V^n_{\omega_n}(t,p)-e^{-\bar A (t-p) })[-\alpha_0^n(\omega_n
p)v_n(p)\\+\hat\varphi_n(\omega_n p,v_n(p))+\gamma_n(\omega_n p)]
\|_{H^1}\, dp,
\endmultline
$$
$$
\multline
I_2(t):= \int_s^t \|e^{- \bar A (t-p)}[-\alpha_0^n(\omega_n p)(v(p)-v_n(p))\\+\hat
\varphi_n(\omega_n p, v(p))- \hat \varphi_n(\omega_n p,v_n(p))]
\|_{H^1}\,dp,
\endmultline
$$
$$
\multline
I_3(t):= \|\int_s^t e^{-\bar A (t-p)}[(\bar a_0 -\alpha_0^n(\omega_n p)) v(p) \\+\Hat{\Bar
f}(v(p)) - \hat \varphi_n (\omega_n p, v(p)) + \gamma_n(\omega_n p)- \bar g]\, dp
\|_{H^1}.
\endmultline
$$
First of all, we have
$$
\multline
\| V^n_{\omega_n}(t,s)u_n -e^{- \bar
A(t-s)}u\|_{H^1} \leq M (t-s)^{-1/2}\|u_n-u\|_{L^2}\\+ \sup_{t \in
[s,s+T]}\|(V^n_{\omega_n}(t,s)-e^{\bar A (t-s)})u
\|_{H^1}.
\endmultline
$$
As for $I_1(t)$, by Proposition \rft linaver..
$$
\multline
I_1(t) \\ \leq \int_s^t (1 +(t-p)^{-1/2})\theta(\omega_n(t-p))\|-\alpha_0^n(\omega_n
p)v_n(p)+\hat\varphi_n(\omega_n p,v_n(p))+\gamma_n(\omega_n p)\|_{L^2}\,dp \\\leq (C \tilde
K+\tilde C(\tilde K + {\tilde K}^{\beta +1})) \int_s^t (t-p)^{-1/2}\theta(\omega_n (t-p))\, dp.
\endmultline
$$
By the same argument used in the proof of Lemma \rft tecnico1.., we find that $I_1(t)\rightarrow
0$ as $n \rightarrow \infty$ uniformly on $[s,s+T]$.\par
Next we consider $I_2(t)$:
$$
\multline
I_2(t)\leq \int_s^t M(1+(t-p)^{-1/2})\|-\alpha_0^n(\omega_n p)(v(p)-v_n(p))\\+\hat
\varphi_n(\omega_n p, v(p))- \hat \varphi_n(\omega_n p,v_n(p))\|_{L^2} \, dp \\ \leq \int_s^t
M(1+(t-p)^{-1/2}) (C\|v(p)-v_n(p)\|_{H^1}+\tilde C(1+2\tilde K^{\beta})\|v(p)-v_n(p)\|_{H^1})\,
dp \\ \leq Q \int_s^t (t-p)^{-1/2}\|v_n(p)-v(p)\|_{H^1}\,dp,
\endmultline
$$
where $Q$ is a positive constant. \par
Finally,
$$
\multline
I_3(t) \leq \| \int_s^t e^{-\bar A (t-p)}[(\bar a_0 -\alpha_0^n(\omega_n p)) v(p)]\,dp\|_{H^1}
\\+\|\int_s^t e^{-\bar A (t-p)}[\Hat{\Bar
f}(v(p)) - \hat \varphi_n (\omega_n p, v(p))]\,dp\|_{H^1}\\+
\|\int_s^t e^{-\bar A (t-p)}[\gamma_n(\omega_n p)- \bar g]\,dp\|_{H^1},
\endmultline
$$
and, by Lemma \rft tecnico1.., $I_3(t)\rightarrow 0$ as $n\rightarrow \infty$ uniformly on
$[s,s+T]$.\par
Summing up, for $t\in ]s,s+T]$ we get
$$
\|v_n(t)-v(t)\|_{H^1}\leq (t-s)^{-1/2}F_n +Q\int_s^t (t-p)^{-1/2}\|v_n(p)-v(p)\|_{H^1}\,dp,
$$
where $F_n \rightarrow 0$ as $n\rightarrow \infty$.\par
By the singular version of Gronwall's inequality (see \cite{\rfa He.., Th. 7.1.1}),
$$
\|v_n(t)-v(t)\|_{H^1}\leq \tilde Q F_n(t-s)^{-1/2},
$$
where $\tilde Q$ is a positive constant. This implies that $v_n(t)\rightarrow v(t)$ in
$H^1(\R^N)$ uniformly on $[s+\delta,s+T]$ for every $\delta >0$, and completes the proof.
\qed\enddemo

\remark{Remark}
If in Theorem \rft contin2.. we assume that $u_n \rightarrow u$ in $H^1(\R^N)$, by the same
techniques we can show that $v_n(t)\rightarrow v(t)$ in $H^1(\R^N)$ uniformly on $[s,s+T]$.
\endremark
\medskip

The following lemma provides a kind of joint asymptotic compactness of
$\Pi^{\sigma}_{\omega}(t,s)$ with respect to $t$ and $\omega$.
\proclaim{Lemma \dft ascomp2..}
Let $(u_n)_{n\in\N}$ be a bounded sequence in $H^1(\R^N)$, $(\sigma_n)_{n\in\N}$ an
arbitrary sequence in $\Sigma$, $(t_n)_{n\in\N}$ and $(\omega_n)_{n\in\N}$ two sequences of
positive real numbers, $t_n\to+\infty$ and $\omega_n\to+\infty$ as $n\to\infty$. Then there
exists
$u_\infty\in H^1(\R^N)$ such that, up to a subsequence,
$$
\Pi^{\sigma_n}_{\omega_n}(t_n,0)u_n\to u_\infty\quad\text{in $H^1(\R^N)$ as $n\to\infty$.}
$$
\endproclaim
\demo{Proof}
The proof is analogous to that of Proposition \rft ascomp1..: from the boundedness
of $(u_n)_{n\in\N}$ in $H^1(\R^N)$ it follows that there exists $u_\infty\in H^1(\R^N)$ such
that, up to a subsequence,
$$
\Pi^{\sigma_n}_{\omega_n}(t_n-1,0)u_n\weakto \cl u_\infty\quad\text{in $H^1(\R^N)$ as
$n\to\infty$.}
$$
Again by Lemma \rft ball1..
$$
\Pi^{\sigma_n}_{\omega_n}(t_n-1,0)u_n\to \cl u_\infty\quad\text{in $L^2(\R^N)$ as
$n\to\infty$,}
$$
and by Theorem \rft contin2.. 
$$
\multline
\Pi^{\sigma_n}_{\omega_n}(t_n,0)u_n
=\Pi^{\sigma_n}_{\omega_n}(t_n,t_n-1)\Pi^{\sigma_n}_{\omega_n}(t_n-1,0)u_n\\
=\Pi^{T_{\omega_n}(t_n-1)\sigma_n}_{\omega_n}(1,0)\Pi^{\sigma_n}_{\omega_n}(t_n-1,0)u_n
\to \pi(1)\cl u_\infty\quad\text{in $H^1(\R^N)$ as $n\to\infty$.}
\endmultline
$$
This completes the proof.
\qed\enddemo

Finally, we can prove the upper-semicontinuity result announced in the Introduction:
\proclaim{Theorem \dft upsem..}
For $\omega >0$, let ${\Cal A}_\omega^\Sigma$ be the $\Sigma$-uniform attractor of the
family of processes $\{\,\Pi^\sigma_\omega\mid \sigma\in\Sigma\,\}$.
Moreover, let ${\Cal A}$ be the attractor of the semiflow $\pi$.
Then for every $\delta >0$
there exists
$\cl
\omega>0$ such that if
$\omega\geq\cl\omega$
$$ \d_{H^1}({\Cal A}^\Sigma_\omega,{\Cal A}):= \max_{u \in {\Cal
A}^\Sigma_\omega}\d_{H^1}(u,{\Cal A})< \delta.$$
\endproclaim

\demo{Proof} Let's assume, by contradiction, that the thesis is not true: then there exist
$\cl
\delta >0$, a sequence
$(\omega_n)_{n \in \N}$
 of positive numbers, $\omega_n \to +\infty$, and a sequence $(u_n)_{n \in \N}$, $u_n \in {\Cal
A}^\Sigma_{\omega_n}$ for all $n\in\N$, such that
$$
\d_{H^1}(u_n, {\Cal A})\geq \cl \delta\quad\text{for all $n\in\N$.}
$$
Since $u_n \in {\Cal A}^\Sigma_{\omega_n}$, by Theorem \rft ChepVish..
 for every $n \in \N$ there exists
$\sigma_n\in\Sigma$ and $v_n\in{\Cal K}^{\sigma_n}_{\omega_n}$ such that $u_n=v_n(0)$.
Since
$$
\left(t\mapsto v_n(t+h)\right)\in{\Cal K}_{\omega_n}^{T_{\omega_n}(h)\sigma_n}$$ for all
$h\in\R$, it follows that
$v_n(t)\in{\Cal A}^{\Sigma}_{\omega_n}$ for all $t\in\R$. Hence,
by Proposition \rft absorb1.., there exists $K>0$, independent of $n$, such that
$\|v_n(t)\|_{H^1}\leq K$ for $t\in \R$.

Let $k$ be a positive integer and let $(h_n)_{n \in \N}$ be a sequence of positive
numbers, $h_n \to +\infty$. We have
$$
v_n(-k)=\Pi^{\sigma_n}_{\omega_n}(-k,-h_n-k)v(-h_n-k)=
\Pi^{T_{\omega_n}(-k-h_n)\sigma_n}_{\omega_n}(h_n,0)v(-h_n-k),
$$
so, by Proposition \rft ascomp2.., there
exists $\cl{u}_k \in H^1(\R^N)$ such that, up to a subsequence,
$$
v_n(-k)\to\cl u_k \quad\text{in $H^1(\R^N)$ as $n \to +\infty$.}\tag\dff sgrunt1..
$$
By a  Cantor diagonal procedure, we can assume that \rff sgrunt1.. holds for
positive integer $k$.
By Theorem \rft contin2.., for every $t>-k$
$$
\multline
v_n(t)=\Pi^{\sigma_n}_{\omega_n}(t,-k)v_n(-k)
=\Pi^{T_{\omega_n}(-k)\sigma_n}_{\omega_n}(t+k,0)v_n(-k)\\
\to \pi(t+k)\cl u_k
\quad\text{in $H^1(\R^N)$ as $n\to\infty$.}\endmultline\tag\dff xyz..
$$
In particular,
choosing $t=0$ we get
$$
u_n\to\pi(k) \cl u_k \quad\text{in $H^1(\R^N)$ as $n\to\infty$.}
$$
Notice that $\pi(k)\cl{u}_k$ is independent of $k$, so we can define
 $u_\infty:= \pi(k)\cl{u}_k$. The proof will be complete if we show that $u_\infty \in {\Cal
A}$. So we must prove that there exists a full bounded solution $v_\infty(t)$ of the
semiflow $\pi$, such that $v_\infty(0)=u_\infty$
To this end, we just have to define $v_\infty(t):= \pi(t+k) \cl{u}_k$, $t> -k$.
By \rff xyz.. it follows that $\pi(t+k)\cl{u}_k$ is independent of $k$ and therefore
$v_\infty(t)$ is unambiguously defined for every $t \in \R$. Moreover, $v_\infty(t)$ is by
construction a full bounded solution of $\pi$, with $v_\infty(0)=u_\infty$. This finally
implies that
$u_\infty\in {\Cal A}$, a contradiction.
\qed
\enddemo

\enddocument\bye